\documentclass[12pt]{article}
\usepackage[ansinew]{inputenc}
\usepackage[spanish,english]{babel}
\usepackage{amsmath,amsthm,amsopn,amstext,amscd,amsfonts,amssymb}
\usepackage{geometry,longtable}
\usepackage{latexsym,enumerate,bm}
\usepackage{graphicx}

\newcommand{\beq}{\begin{equation}}
\newcommand{\eeq}{\end{equation}}
\newcommand{\bed}{\begin{displaymath}}
\newcommand{\eed}{\end{displaymath}}

\def \bbeta {\boldsymbol\beta}

\def \btheta {\boldsymbol \theta}
\def \eeta {\boldsymbol \eta}
\def \bepsilon {\boldsymbol \epsilon}

\def \bDelta {\boldsymbol\Delta}
\def \bkappa {\boldsymbol\kappa}
\def \bmu {\boldsymbol\mu}
\def \bnu {\boldsymbol\nu}

\def \bvarsigma {\boldsymbol\varsigma}
\def \ds {\displaystyle}

\def \cB {{\cal B}}

\def \cN {{\cal N}}
\def \cF {{\cal F}}

\def \RR {I\!\!R}

\def \bA {{\mathbf A}}
\def \ba {{\mathbf a}}
\def \bB {{\mathbf B}}
\def \bb {{\mathbf b}}
\def \bc {{\mathbf c}}

\def \bd {{\mathbf d}}
\def \bE {{\mathbf E}}
\def \be {{\mathbf e}}
\def \ff {{\mathbf f}}
\def \bF {{\mathbf F}}
\def \bg {{\mathbf g}}
\def \bG {{\mathbf G}}

\def \bh {{\mathbf h}}
\def \bI {{\mathbf I}}

\def \br {{\mathbf r}}
\def \bs {{\mathbf s}}

\def \bP {{\mathbf P}}
\def \bq {{\mathbf q}}
\def \bQ {{\mathbf Q}}

\def \bm {{\mathbf m}}
\def \bX {{\mathbf X}}
\def \bx {{\mathbf x}}
\def \bu {{\mathbf u}}
\def \bv {{\mathbf v}}
\def \bV {{\mathbf V}}
\def \by {{\mathbf y}}
\def \bZ {{\mathbf Z}}
\def \bz {{\mathbf z}}
\def \bW {{\mathbf W}}

\def \diag {\mbox{diag}}
\def \tr {\mbox{tr}}
\def \col {\mbox{col}}
\def \uno {\mathbf 1}
\def \cero {\mathbf 0}

\newtheorem{theo}{Theorem}
\newtheorem{coro}{Corollary}
\newtheorem{lemm}{Lemma}
\newtheorem{prop}{Proposition}

\newtheorem{rema}{Remark}

\title{\bf Best prediction under a nested error model with log transformation\footnote{Supported by the Spanish grants SEJ2007-64500, MTM2012-37077-C02-01, MTM-2012-33740 and ECO-2011-25706. Part of this work was done during a research stay of the second author in the Institute of Statistics of the University of Neuchâtel.} }
\author{Nirian Martín
\ and Isabel Molina
\\
Department of Statistics, Universidad Carlos III de Madrid}
\date{}

\begin{document}
\maketitle

\begin{quote}
\vspace{0.3 cm}\\
{\bf Abstract:} In regression models involving economic variables
such as income, log transformation is typically taken to
achieve approximate normality and stabilize the variance. However,
often the interest is predicting individual values or means
of the variable in the original scale. Back transformation of predicted values
introduces a non-negligible bias. Moreover, assessing the
uncertainty of the actual predictor is not straightforward. In
this paper, a nested error model for the log transformation of the
target variable is considered. Nested error models are widely used for
estimation of means in subpopulations with small sample sizes (small areas),
by linking all the areas through common parameters. These common parameters are estimated
using the overall set of sample data, which leads to much more efficient small area
estimators. Analytical expressions for the best predictors of individual values of the original variable and of
small area means are obtained under the nested
error model with log transformation of the target variable. Empirical
best predictors are defined by estimating the unknown model
parameters in the best predictors. Exact mean squared errors of the best predictors and
second order approximations to the mean squared errors of the
empirical best predictors are derived. Mean
squared error estimators that are second order correct are also obtained.
An example with Mexican data
on living conditions illustrates the procedures.
\end{quote}

\section{Introduction}

In Econometric regression models, variables such as income or expenditure
are often transformed with logarithm to achieve homoscedastic errors with
approximately normal distribution. However, the variable of interest remains
to be the untransformed one. Target characteristics of the study variable
such as the values for out-of-sample individuals or the means for specific
subpopulations become then functions of the exponentials of the dependent
variable in the model. However, the predictors obtained by transforming back
the individual predicted values are biased. Usual bias-corrections are only
approximations and optimality properties are lost. However, the exact
expression for the optimal predictors can be obtained analytically for
certain models. A model that is often used for small area estimation is the
nested-error linear regression model proposed by Battesse, Harter and Fuller
(1988) to estimate the area under production of corn and soybeans in a
number of counties. In small area estimation, the lack of sample
observations in some of the areas of interest is solved by linking all areas
through the common regression parameters but including at the same time
random area effects that represent the unexplained between area variation.
The common parameters are estimated using the sample observations from all
the areas together and this leads to great efficiency gains with respect to
estimators that use only the area-specific sample data (direct estimators).
This kind of model is used in Econometric applications as well, see e.g.
Elbers, Lanjouw and Lanjouw (2003) or Molina and Rao (2010), who employed
this model to estimate poverty indicators in small areas. For more details
on small area estimation methods, see the monograph by Rao and Molina (2015) and the
recent review by Pfeffermann (2013).

Assessing the reliability, or uncertainty, of the obtained predictors is
crucial in practical applications. A popular uncertainty measure is the mean
squared error (MSE), also called mean squared prediction error. MSEs of
optimal predictors of small area parameters have been obtained under certain
models but only for simple parameters, see e.g. Das, Jiang and Rao (2004).
The MSE of an individual prediction under a nested-error model with
log-transformation that is second-order correct has not been obtained yet. Moreover, when predicting the
mean of the original variable in a given area, the optimal predictor is
function of the predicted values for the out-of-sample individuals from that
area. Since the individuals belong to the same area, due to the presence of
the area effects, individual predictors are not independent. Then mean
crossed product errors (MCPEs) between pairs of individual predictions are
needed to derive the MSE of the predictor of the mean in that area.

Here we obtain optimal predictors for individual values of the target
variable in out-of-sample units and also for small area means. Additionally,
second-order asymptotic approximations for the MCPEs of pairs of individual
predictions are derived, which lead to good approximations for the MSEs of
predicted area means. In the small area estimation literature, this was done
previously only under area-level models by Slud and Maiti (2006). Under a
unit-level model, Molina (2009) dealt with estimation of exponentials of
mixed effects, i.e. exponentials of linear functions of the fixed and the
random effects in the model; the individual values of the original variable
cannot be expressed as special cases of these parameters. Thus, the target
parameters and not the same and consequently results are also different. In
particular, certain crossed-product terms appearing in the MCPE that are of
lower order in Molina (2009), are not negligible when predicting individual
observations. In fact, those crossed-product terms are typically neglected
in small area estimation applications. Here we show that these terms cannot
be neglected and give their analytical expression up to $o(D^{-1})$ terms,
where $D$ is the number of areas.

Analytical approximations for the uncertainty measures have a complex shape
and users might prefer to use resampling procedures such as bootstrap
methods. Gonz\'alez-Manteiga et al. (2008) proposed a parametric bootstrap
method designed for finite populations under a nested error model that is
suitable in this paper. However, Gonz\'alez-Manteiga et al. (2008) proved
consistency of the bootstrap MSE estimator when the target parameters are
linear. For our particular non-linear parameters, consistency remains to be
proved. Nevertheless, once an analytical asymptotic expression is available
for the true MSE, the technique of imitation used in that paper can be
followed to achieve the consistency in this paper. Thus, the theoretical
results for the MSE approximation that are obtained in this paper lead
automatically to the consistency of the corresponding bootstrap MSE
estimators.

The paper is organized as follows. The considered model and the target
quantities are introduced in Section \ref{secModel}. This section also gives the best predictor and
first and second-stage empirical best predictors of the target quantities. Section \ref{fitting} describes
usual likelihood-based fitting methods.
MCPEs and MSEs of first-stage empirical best predictors are
obtained in Section \ref{secUncertainty}, and for second-stage empirical
best predictors, second-order approximations to the analogous uncertainty
measures are given in Section \ref{MSEEB2}. Second-order unbiased estimators
of these uncertainty measures are provided in Section \ref{secEstimUncer}.
Section \ref{SecBoot} describes a parametric bootstrap procedure for
estimation of the uncertainty. Section \ref{sec:simexp} describes the result of
a simulation experiment comparing the proposed predictor with existing ones.
Section \ref{secExample} illustrates the
procedures through the estimation of mean income in municipalities from Mexico.
Finally, the proofs of all the theorems are included in the Appendix.

\section{Model, target quantities and predictors}

\label{secModel}

When estimating characteristics of subpopulations that have varying sizes,
it seems convenient to work under a finite population setup. Here we
consider that the population $U$ is finite and contains $N$ units. This
population is partitioned into $D$ subpopulations $U_1,\ldots,U_D$, also
called areas or domains, of sizes $N_1,\ldots,N_D$. The data is obtained
from a sample $s$ of size $n$ drawn from the population $U$. We denote by $%
s_d$ the subsample from domain $d$, of (fixed) size $n_d$, $d=1,\ldots,D$,
with $\sum_{d=1}^Dn_d=n$, and by $\bar s_d=U_d-s_d$ the sample complement
from area $d$, of size $N_d-n_d$, $d=1,\ldots,D$.

The goal is to predict the value $w_{di}$ of the variable of interest for an
out-of-sample individual $i$ within area $d$, or the area mean $%
N_{d}^{-1}\sum_{i=1}^{N_{d}}w_{di}$, based on a regression model for $w_{di}$%
. If $w_{di}$ represents a measurement of an economical variable such as
income or expenditure, it is customary to consider the logarithm of $w_{di}$
as dependent variable in a regression model. Moreover, in many applications,
the available auxiliary variables do not explain sufficiently well all the
between-area variation that data exhibit. Then, random area effects
representing this unexplained variation are included in the model. This is
typically done in small area estimation applications. Here we assume the
following linear regression model with random area effects, also known as
nested-error model, for the log-transformed variables $y_{di}=\log w_{di}$,
\begin{equation} \label{neemodel}
y_{di}=\bx_{di}^{\prime }\bbeta+u_{d}+e_{di},\ u_{d}\overset{iid}{\sim }\cN%
(0,\sigma _{u}^{2}),\  e_{di}\overset{iid}{\sim }\cN(0,\sigma _{e}^{2}),\quad i=1,\ldots ,N_{d},\
d=1,\ldots ,D.
\end{equation}%
Here, $\bx_{di}$ is a vector containing the values of $p$ explanatory
variables for $i$-th individual in $d$-th area, $\bbeta\in \RR^{p}$ is the
vector of unknown regression coefficients, $e_{di}$ is the individual error,
$u_{d}$ is the random effect of area $d$, with random effects $\{u_{d}\}$
and errors $\{e_{di}\}$ assumed to be independent, and finally $\sigma
_{u}^{2}$ and $\sigma _{e}^{2}$ are the unknown random effects and
individual error variances respectively, called variance components. We
denote by $\btheta=(\sigma _{u}^{2},\sigma _{e}^{2})^{\prime }$ the vector
of variance components and by $\Theta =\{(\sigma _{u}^{2},\sigma
_{e}^{2})^{\prime };\sigma _{u}^{2}\geq 0,\ \sigma _{e}^{2}>0\}$ the space
where these parameters lie. Notation $\bbeta$ and $\btheta$ will refer
hereafter to generic elements from $\RR^{p}$ and $\Theta $, whereas $\bbeta%
_{0}$ and $\btheta_{0}$ will be the respective true values of $\bbeta$ and $%
\btheta$, where $\btheta_{0}$ is supposed to be within the interior of $%
\Theta $. For a quantity $A(\bbeta,\btheta)$ depending on $\btheta$ and/or $%
\bbeta$, we will use many times the notation $A$, omitting the explicit
dependence on $\bbeta$ and/or $\btheta$.

If we intend to estimate the mean of an area with a poor sample size $n_d$,
the estimators that use only the $n_d$ area-specific observations, called
direct estimators, are highly inefficient. Model (\ref{neemodel}) links all
the areas through the common parameters $\bbeta$, $\sigma_u^2$ and $%
\sigma_e^2$, which allows us to ``borrow strength" from all the areas when
estimating a particular area mean. However, even though the model is assumed
for $y_{di}=\log w_{di}$, the target parameter remains to be the area mean
of the untransformed variables, which can be expressed in terms of the
dependent variables in the model as
\begin{equation*}
\tau_d=\frac{1}{N_d}\sum_{i=1}^{N_d}w_{di}=\frac{1}{N_d}\sum_{i=1}^{N_d}%
\exp(y_{di}),\quad d=1,\ldots,D.
\end{equation*}





Here we intend to estimate single values $w_{di}=\exp(y_{di})$ of the target variable in
out-of-sample units $i\in \bar s_d$ and area means $\tau_d=N_d^{-1}\sum_{i=1}^{N_d}\exp(y_{di})$, when the variables $y_{di}$ in the population
units follow model (\ref{neemodel}). These target quantities are special cases of a
general parameter of the form $h(\by_{d})$, where $h(\cdot)$ is a measurable function and
$\by_{d}=(y_{d1},\ldots ,y_{dN_{d}})^{\prime }$ is the vector of outcomes for domain $d$. Defining also
$\bX_{d}=(\bx_{d1},\ldots,\bx_{dN_{d}})^{\prime }$ and $\be_{d}=(e_{d1},\ldots ,e_{dN_{d}})^{\prime }$,
the model reads
\begin{equation}\label{neemodelmat}
\by_{d}=\bX_{d}\bbeta+u_{d}\uno_{N_{d}}+\be_{d},\ u_{d}\overset{iid}{\sim }%
\cN(0,\sigma _{u}^{2}), \ \be_{d}\overset{ind}{\sim }\cN_{N_{d}}(\cero_{N_{d}},\sigma _{e}^{2}\bI%
_{N_{d}}),\quad d=1,\ldots ,D,
\end{equation}%
where $\cero_{k}$ is a $k$-vector of zeros, $\uno_{k}$
is a $k$-vector of ones and $\bI_{k}$ is the $k\times k$ identity
matrix. The covariance matrix of $\by_{d}$ is equal to $\bV_{d}=\sigma
_{u}^{2}\uno_{N_{d}}\uno_{N_{d}}^{\prime }+\sigma _{e}^{2}\bI_{N_{d}}=\bV%
_{d}(\btheta)$. Let us arrange the elements from domain $d$ into sample and out-of-sample elements, as
\begin{equation*}
\by_{d}=\left(
\begin{array}{c}
\by_{ds} \\
\by_{dr}%
\end{array}%
\right) ,\quad \bX_{d}=\left(
\begin{array}{c}
\bX_{ds} \\
\bX_{dr}%
\end{array}%
\right) ,\quad \bV_{d}=\left(
\begin{array}{cc}
\bV_{ds} & \bV_{dsr} \\
\bV_{drs} & \bV_{dr}%
\end{array}%
\right) .
\end{equation*}%
The ``best predictor" $\tilde\delta_d$ of a general parameter $\delta_d=h(\by_d)$ is the function
of the sample data $\by_{ds}$ with minimum mean squared error $\mbox{MSE}(\tilde\delta_d)=E(\tilde\delta_d-\delta_d)^2$ and
is given by $\tilde\delta_d=E_{\by_{dr}}\{h(\by_{d})|\by_{ds}\}$, where the expectation is taken with respect to the distribution of $\by_{dr}|\by_{ds}$.
The best predictor is exactly unbiased in the sense $E_{\by_{ds}}(\tilde\delta_d)=E_{\by_d}(\delta_d)$.
Since by (\ref{neemodelmat}) we have $\by_{d}\sim \cN(\bX_{d}\bbeta,\bV_{d})$, the desired conditional distribution is
\begin{equation}
\by_{dr}|\by_{ds}\overset{ind}{\sim }\cN_{N_{d}-n_{d}}(\bmu_{dr|s},\bV%
_{dr|s}),\quad d=1,\ldots ,D,  \label{condist}
\end{equation}%
with mean vector and covariance matrix given by
\begin{equation*}
\bmu_{dr|s} =\bX_{dr}\bbeta+\bV_{drs}\bV_{ds}^{-1}(\by_{ds}-\bX_{ds}\bbeta),
\quad
\bV_{dr|s} =\bV_{dr}-\bV_{drs}\bV_{ds}^{-1}\bV_{dsr}.
\end{equation*}%
Under the nested-error model (\ref{neemodel}), they reduce to
\begin{align}
\bmu_{dr|s}& =\bX_{dr}\bbeta+\uno_{N_{d}-n_{d}}\gamma _{d}(\bar{y}_{ds}-\bar{%
\bx}_{ds}^{\prime }\bbeta)  \label{murs} \\
\bV_{dr|s}& =\sigma _{u}^{2}(1-\gamma _{d})\uno_{N_{d}-n_{d}}\uno%
_{N_{d}-n_{d}}^{\prime }+\sigma _{e}^{2}\bI_{N_{d}-n_{d}},  \label{vrs}
\end{align}%
where $\bar{y}_{ds}=n_{d}^{-1}\sum_{i\in s_{d}}y_{di}$, $\bar{\bx}_{ds}=
n_{d}^{-1}\sum_{i\in s_{d}}\bx_{di}$ and $\gamma _{d}=\sigma_{u}^{2}/(\sigma _{u}^{2}+\sigma _{e}^{2}/n_{d})$.

Based on the conditional distribution (\ref{condist}) with mean vector given
in (\ref{murs}) and covariance matrix (\ref{vrs}), the next theorem gives
closed-form expressions for the best predictors of $w_{di}=\exp(y_{di})$ and
$\tau_d=N_d^{-1}\sum_{i=1}^{N_d}\exp(y_{di})$.

\begin{theo}
\label{expresBP} Under the nested-error model with log-transformation (\ref%
{neemodel}), it holds:

\begin{itemize}
\item[(i)] The best predictor of $w_{di}=\exp(y_{di})$, for $i\in \bar s_d$,
is given by
\begin{equation}  \label{bpw}
\tilde w_{di}=\tilde w_{di}(\bbeta,\btheta)=\exp(\tilde y_{di}+\alpha_d),
\end{equation}
where
$\tilde y_{di}=\bx_{di}^{\prime }\bbeta+\gamma_d(\bar y_{ds}-\bar \bx
_{ds}^{\prime }\bbeta)$ and $\alpha_d=\{\sigma_u^2(1-
\gamma_d)+\sigma_e^2\}/2$.
\vspace{0.2 cm}
\item[(ii)] The best predictor of $\tau_d=N_d^{-1}\sum_{i=1}^{N_d}%
\exp(y_{di})$ is given by
\begin{equation}  \label{BPtaud}
\tilde \tau_d=\tilde \tau_d(\bbeta,\btheta)=\frac{1}{N_d}\left(\sum_{i\in
s_d}w_{di}+\sum_{i\in \bar s_d}\tilde w_{di}\right).
\end{equation}
\end{itemize}
\end{theo}

\begin{rema}
\label{nonsampleneed} In contrast with the case of estimation of a small
area mean under a nested error model without log-transformation, the best
predictor of the small area mean $\tau_d$ given in (\ref{BPtaud}) requires
the values of the auxiliary variables $\bx_{di}$ for each out-of-sample unit
$i\in \bar s_d$ and not only of area totals or means of the auxiliary
variables. Censuses of potentially useful auxiliary variables are available for practically all European countries and many other countries all over the world.
\end{rema}

Molina (2009) proposed the bias-corrected predictor $\tilde w_{di}^M=\exp(\tilde y_{di}+\alpha_d^M)$, where $\alpha_d^M= \sigma_u^2(1-\gamma_d)/2$, which
is similar to the best predictor $\tilde w_{di}$ given in (\ref{bpw}). However, they are not exactly the same because the target
parameters in Molina (2009) are of the type $\exp(\bx_{di}^{\prime }\bbeta+u_d)$, which differ from our target parameters here given by the individual observations $w_{di}=\exp(y_{di})=\exp(\bx_{di}^{\prime }\bbeta+u_d+e_{di})$.
Nevertheless, it is interesting to study how Molina (2009)'s predictor $\tilde w_{di}^M$ performs for $w_{di}=\exp(y_{di})$. The next result gives the relative bias of $\tilde w_{di}^M$ and of the naive predictor obtained by back-transforming the predicted model responses, $\tilde w_{di}^N=\exp(\tilde y_{di})$. By this result, these two predictors are negatively biased unlike the best predictor $\tilde w_{di}$ given in (\ref{bpw}) and $|RB(\tilde w_{di}^N)|\geq |RB(\tilde w_{di}^M)|$.

\begin{prop} \label{RBnaiveM} Under model \eqref{neemodel}, it holds:
\begin{itemize}
\item[(i)] $RB(\tilde w_{di}^N)=\exp(\alpha_d)-1$;
\item[(ii)] $RB(\tilde w_{di}^M)=\exp(\sigma_e^2/2)-1$.
\end{itemize}
\end{prop}

The best predictors $\tilde{w}_{di}(\bbeta,\btheta)$ and
$\tilde{\tau}_{d}(\bbeta,\btheta)$ depend on the true values of
$\bbeta$ and $\btheta$, which are unknown in practice. Next we
define first and second-stage empirical best (EB) predictors obtained by
estimating these unknown parameters in two stages. First, define
the following vectors and matrices containing the sample elements from all
the areas
\begin{align*}
& \by_{s}=(\by_{1s}^{\prime },\ldots ,\by_{Ds}^{\prime })^{\prime },\quad \bX%
_{s}=(\bX_{1s}^{\prime },\ldots ,\bX_{Ds}^{\prime })^{\prime },\quad \be%
_{s}=(\be_{1s}^{\prime },\ldots ,\be_{Ds}^{\prime })^{\prime }, \\
& \bZ_{s}=\diag_{1\leq d\leq D}(\uno_{n_{d}}),\quad \bu=(u_{1},\ldots
,u_{D})^{\prime }.
\end{align*}%
Then, the model for the sample units can be written as
\begin{equation*}
\by_{s}=\bX_{s}\bbeta+\bZ_{s}\bu+\be_{s},\quad \bu\sim \cN_{D}(\cero%
_{D},\sigma _{u}^{2}\bI_{D}),\quad \be_{s}\sim \cN_{n}(\cero_{n},\sigma
_{e}^{2}\bI_{n}),
\end{equation*}%
and the covariance matrix of $\by_{s}$ is given by
\begin{equation*}
\bV_{s}=\diag_{1\leq d\leq D}(\bV_{ds}),\quad \bV_{ds}=\sigma _{u}^{2}\uno%
_{n_{d}}\uno_{n_{d}}^{\prime }+\sigma _{e}^{2}\bI_{n_{d}}.
\end{equation*}

The first-stage EB predictor is obtained under the assumption that $\btheta$
is known but $\bbeta$ is unknown. The maximum likelihood (ML) estimator of $%
\bbeta$ under normality, which is also the weighted least squares (WLS)
estimator of $\bbeta$ without normality reads
\begin{equation}  \label{betaWLS}
\tilde\bbeta(\btheta)=(\bX_s^{\prime }\bV_s^{-1}\bX_s)^{-1}\bX_s^{\prime }\bV%
_s^{-1}\by_s.
\end{equation}
The first-stage EB predictors of $w_{di}$ and $\tau_d$ are then
\begin{equation}  \label{EB1}
\hat w_{di}=\hat w_{di}(\btheta)=\tilde w_{di}(\tilde\bbeta(\btheta),\btheta%
),\quad \hat\tau_d=\hat\tau_d(\btheta)=\tilde\tau_d(\tilde\bbeta(\btheta),%
\btheta).
\end{equation}
The next result gives asymptotic unbiasedness of $\hat w_{di}$ at the log scale.
\begin{prop}\label{proplog}
Under model \eqref{neemodel} and assumptions (H1)--(H3) of Section \ref{MSEEB2}, it holds
$$
\log E(\hat w_{di})=\log E(w_{di})+O(D^{-1}).
$$
\end{prop}

Finally, the second-stage EB predictors of $w_{di}$ and $\tau_d$ are
obtained by replacing the unknown $\btheta$ in (\ref{EB1}) by a consistent
estimator $\hat\btheta$, that is,
\begin{equation}  \label{EB2}
\hat w_{di}^E=\hat w_{di}(\hat\btheta)=\tilde w_{di}(\tilde\bbeta(\hat\btheta%
),\hat\btheta),\quad \hat\tau_d^E=\hat\tau_d(\hat\btheta)=\tilde\tau_d(\tilde%
\bbeta(\hat\btheta),\hat\btheta).
\end{equation}
Section \ref{MSEEB2} describes typical methods for consistent estimation of $\btheta$ under model \eqref{neemodel}.

\section{Fitting methods}\label{fitting}

 A typical estimation method is maximum likelihood (ML), which provides consistent and asymptotically efficient estimators of the variance components
 (Miller, 1973). The ML estimator $\hat\btheta=(\hat\sigma_u^2,\hat\sigma_e^2)'$ of $\btheta=(\sigma_u^2,\sigma_e^2)'$ maximizes the penalized log-likelihood, given by
 \begin{equation}\label{pllike}
 l_P(\btheta)=c-\frac{1}{2}\left(\log|\bV_s|+\by_s'\bP_s\by_s\right),\quad
 \bP_s=\bV_s^{-1}-\bV_s^{-1}\bX_s\bQ_s\bX_s'\bV_s^{-1},
 \end{equation}
 where $c$ denotes a generic constant. The score vector is defined as $\bs(\btheta)=\partial l_P(\btheta)/\partial \btheta=(s_1(\btheta),s_2(\btheta))'$. In terms of the vector $\bv_s=\by_s-\bX_s\bbeta=\bZ_s\bu+\be_s$,
 the elements of the score vector are given by
 \beq\label{score}
 s_h(\btheta)=-\frac{1}{2}\tr(\bV_s^{-1}\bDelta_h)+\frac{1}{2}\bv_s^{\prime}\bP_s\bDelta_h\bP_s\bv_s,\quad h=1,2,
 \eeq
 where $\bDelta_h=\partial \bV_s/\partial \theta_h$, that is, $\bDelta_1=\bZ_s\bZ_s^{\prime}$ and $\bDelta_2=\bI_n$.
 The ML estimator of $\btheta$ is then obtained solving the equation system $\bs(\btheta)=\cero_2$ together with equation \eqref{betaWLS} for $\bbeta$. Since equations are non-linear, numerical algorithms such as Newton-Raphson or Fisher-Scoring are typically applied. These algorithms require respectively the elements of Hessian matrix or the Fisher information matrix. The Hessian matrix is defined as $H(\btheta)=\partial^2 l_P(\btheta)/\partial \btheta^2=(H_{h\ell}(\btheta))$, where
 $$
 H_{h\ell}(\btheta)=\frac{1}{2}\tr(\bV_s^{-1}\bDelta_h\bV_s^{-1}\bDelta_{\ell})-\bv_s'\bP_s\bDelta_h\bP_s\bDelta_{\ell}\bP_s\bv_s,\quad h,\ell=1,2.
 $$
 Finally, the Fisher information matrix is $\cF(\btheta)=E\{-H(\btheta)\}=(\cF_{h\ell}(\btheta))$, where
 $$
 \cF_{h\ell}(\btheta)=-\frac{1}{2}\tr(\bV_s^{-1}\bDelta_h\bV_s^{-1}\bDelta_{\ell})+\tr(\bP_s\bDelta_h\bP_s\bDelta_{\ell}),\quad h,\ell=1,2.
 $$

 A drawback of ML estimator of $\btheta$ is that is does not account for the degrees of freedom due to estimation of $\bbeta$. Restricted ML (REML) corrects for this problem, providing estimators with bias of lower order. This is achieved by transforming the data $\by$ as $\bF'\by$, where $\bF$ is any $n\times (n-p)$ matrix with rank $n-p$ and satisfying $\bF'\bX=\cero_{n-p}$. The REML estimator is the value of $\btheta$ maximizing the so called restricted log-likelihood $l_R$, which is the logarithm of the joint density function of the transformed data $\bF'\by$. Noting that $\bF(\bF'\bV_s\bF)^{-1}\bF'=\bP_s$ (Searle et al. 1992, p.451), this function can be written as
 \begin{equation}\label{pllikeR}
 l_R(\btheta)=c-\frac{1}{2}\left(\log|\bF'\bV_s\bF|+\by_s'\bP_s\by_s\right).
 \end{equation}
 The score vector obtained from $l_R$ is $\bs_R(\btheta)=\partial l_R(\btheta)/\partial \btheta=(s_{R,1}(\btheta),s_{R,2}(\btheta))'$. Using again the relation $\bF(\bF'\bV_s\bF)^{-1}\bF'=\bP_s$, the elements of $\bs_R$ can be expressed as
 \beq\label{score}
 s_{R,h}(\btheta)=-\frac{1}{2}\tr(\bP_s^{-1}\bDelta_h)+\frac{1}{2}\bv_s'\bP_s\bDelta_h\bP_s\bv_s,\quad h=1,2.
 \eeq
 The Hessian matrix obtained from $l_R$ is $H_R(\btheta)=\partial^2 l_R(\btheta)/\partial \btheta^2=(H_{R,h\ell}(\btheta))$, where
 $$
 H_{R,h\ell}(\btheta)=\frac{1}{2}\tr(\bP_s^{-1}\bDelta_h\bP_s^{-1}\bDelta_{\ell})-\bv_s'\bP_s\bDelta_h\bP_s\bDelta_{\ell}\bP_s\bv_s,\quad h,\ell=1,2.
 $$
 Finally, the corresponding Fisher information matrix is in this case given by $\cF_R(\btheta)=E\{-H_R(\btheta)\}=(\cF_{R,h\ell}(\btheta))$, with elements
 $$
 \cF_{R,h\ell}(\btheta)=\frac{1}{2}\tr(\bP_s\bDelta_h\bP_s\bDelta_{\ell}),\quad h,\ell=1,2.
 $$

\section{Uncertainty of first-stage EB predictors}

\label{secUncertainty}

The reliability of a point predictor is typically assessed by its MSE.
When estimating a small area mean $\tau_d$, in virtue of (\ref{BPtaud}),
the MSE of a predictor $\tilde\tau_d$ can be directly
obtained as a function of the MCPEs of pairs of predictors $\hat
w_{di}$ and $\hat w_{dj}$ for out-of-sample units $i,j\in \bar s_d$.
For this reason, in the following we focus on giving the
expressions for the MCPEs of pairs of individual predictors.

Theorem \ref{mseeb} spells out the MCPE of the best predictors $\tilde
w_{di} $ and $\tilde w_{dj}$ for out-of-sample units $i,j\in \bar s_d$,
defined by $\mbox{MCPE}(\tilde w_{di},\tilde w_{dj})=E\{(\tilde
w_{di}-w_{di})(\tilde w_{dj}-w_{dj})\}$. The mean squared error (MSE) of the
best predictor of a single out-of-sample observation $\mbox{MSE}(\tilde
w_{di})=E(\tilde w_{di}-w_{di})^2$, $i\in \bar s_d$ is then obtained taking $%
i=j$. For the area mean $\tau_d$, the MSE of the best predictor $\mbox{MSE}%
(\tilde \tau_d)=E(\tilde \tau_d-\tau_d)^2$ is given in Corollary \ref%
{mseebtau}. Let $1_{\{i=j\}}$ be equal to 1 if $i=j$ and 0 otherwise, and
\begin{equation*}
S_1=\sum_{i\in \bar{s}_{d}}\sum_{j\in \bar{s}_{d},j>i} \exp \left\{ (\bx%
_{di}+\bx_{dj})^{\prime }\bbeta\right\},\quad S_2=\sum_{i\in \bar{s}%
_{d}}\exp \left\{ 2\bx_{di}^{\prime }\bbeta\right\}.
\end{equation*}

\begin{theo}
\label{mseeb} Under the nested-error model with log-transformation (\ref%
{neemodel}), the mean crossed product error of the best predictors $\tilde{w}%
_{di}$ and $\tilde{w}_{dj}$ of $w_{di}$ and $w_{dj}$, for $i,j\in \bar{s}%
_{d} $, is given by
\begin{align*}
\mbox{MCPE}(\tilde{w}_{di},\tilde{w}_{dj})& =\exp \left\{ 2\sigma
_{u}^{2}+\sigma _{e}^{2}+(\bx_{di}+\bx_{dj})^{\prime }\bbeta\right\} \\
& \times \left[1+\left\{ \exp (\sigma _{e}^{2})-1\right\} 1_{\{i=j\}}-\exp \left\{ -\sigma _{u}^{2}(1-\gamma _{d})\right\}
\right] .
\end{align*}
\end{theo}

\begin{coro}
\label{mseebtau} The mean squared error of the best predictor $\tilde{\tau}%
_{d}$ of $\tau _{d}$ is given by
\begin{align*}
\mbox{MSE}(\tilde{\tau}_{d})& =N_{d}^{-2}\exp \left( 2\sigma _{u}^{2}+\sigma
_{e}^{2}\right) \left( 2\left[1- \exp \left\{-\sigma _{u}^{2}(1-\gamma
_{d})\right\}\right] S_{1}\right. \\
& +\left. \left[ \exp(\sigma _{e}^{2})-\exp \left\{ -\sigma _{u}^{2}(1-\gamma _{d})\right\}\right] S_{2}\right) .
\end{align*}
\end{coro}

For a pair of first-stage EB predictors obtained by estimating $\bbeta$
using the WLS estimator given in (\ref{betaWLS}) but assuming that $\btheta$
is known, Theorem \ref{mseeb1} gives the MCPE. The MSE of a single
first-stage EB predictor is obtained setting $j=i$. The following notation
is required:
\begin{equation*}
\bQ_s=(\bX_s^{\prime }\bV_s^{-1}\bX_s)^{-1},\quad h_{d,ij}=\bx_{di}^{\prime }%
\bQ_s\bx_{dj},\quad h_{d,i}=\bx_{di}^{\prime }\bQ_s\bar\bx_{ds},\quad
h_d=\bar\bx_{ds}^{\prime }\bQ_s\bar\bx_{ds}.
\end{equation*}

\begin{theo}
\label{mseeb1} Under the nested-error model with log-transformation (\ref%
{neemodel}), the mean crossed product error of the first-stage EB predictors
$\hat{w}_{di}$ and $\hat{w}_{dj}$, for $i,j\in \bar{s}_{d}$, is given by
\begin{align}
& \mbox{MCPE}(\hat{w}_{di},\hat{w}_{dj})=\exp \left\{ 2\sigma
_{u}^{2}+\sigma _{e}^{2}+(\bx_{di}+\bx_{dj})^{\prime }\bbeta%
\right\}  \label{MCPEEB1} \\
&  \times \left[ 1+\left\{ \exp (\sigma _{e}^{2})-1\right\} 1_{\{i=j\}}+\exp \left\{
(h_{d,ii}+h_{d,jj})/2+h_{d,ij}-2\gamma _{d}^{2}h_{d}-\sigma _{u}^{2}(1-\gamma
_{d})\right\} \right.  \notag \\
&  -\exp \left\{ (h_{d,jj}-\gamma
_{d}^{2}h_{d})/2+\gamma _{d}(h_{d,j}-\gamma _{d}h_d)-\sigma_u^2(1-\gamma_d)\right\}  \notag \\
&  \left. -\exp \left\{ (h_{d,ii}-\gamma
_{d}^{2}h_{d})/2+\gamma_d(h_{d,i}-\gamma _{d}h_d)-\sigma_u^2(1-\gamma_d)\right\} \right] =:M_{1d,ij}(\bbeta,\btheta).  \notag
\end{align}
\end{theo}

\section{Uncertainty of second-stage EB predictors}\label{MSEEB2}

 In practice, the vector of variance components $\btheta=(\sigma_u^2,\sigma_e^2)'$ is also unknown.
 Estimation of $\btheta$ to obtain second-stage EB predictors entails an
 increase in uncertainty and this increase should be accounted for in the MCPE. The additional uncertainty depends on the estimation method
 used for $\btheta$. This section gives an approximation up to $o(D^{-1})$ terms for the MCPE of pairs of individual second-stage EB predictors when model parameters are estimated by ML or REML.

 For the second-stage EB predictors $\hat w_{di}^E=\hat w_{di}(\hat\btheta)$ and $\hat w_{dj}^E=\hat w_{dj}(\hat\btheta)$ of $w_{di}$ and $w_{dj}$, for $i,j\in \bar s_d$, the MCPE can be decomposed as
\begin{align}
& \mbox{MCPE}(\hat{w}_{di}^{E},\hat{w}_{dj}^{E})=\mbox{MCPE}(\hat{w}_{di},%
\hat{w}_{dj})+E\left\{ (\hat{w}_{di}^{E}-\hat{w}_{di})(\hat{w}_{dj}^{E}-\hat{%
w}_{dj})\right\}  \label{mcpeeb2} \\
& \quad +E\left\{ (\hat{w}_{di}^{E}-\hat{w}_{di})(\hat{w}_{dj}-w_{dj})\right%
\} +E\left\{ (\hat{w}_{di}-w_{di})(\hat{w}_{dj}^{E}-\hat{w}_{dj})\right\} .
\notag
\end{align}%
The first term on the right-hand side of (\ref{mcpeeb2}) is already given in
Theorem \ref{mseeb1} above. The remaining terms will be approximated up to $%
o(D^{-1})$ terms under the following assumptions, where $\lambda _{\min }(A)$
denotes the minimum eigenvalue of $A$:

\begin{enumerate}
\item[(H1)] $p<\infty$, $\ds\limsup_{D\to\infty}\max_{1\leq d\leq D}
n_d<\infty$ and $\ds\liminf_{D\to\infty}\min_{1\leq d\leq D} n_d>0$;

\item[(H2)] The elements of the matrix $\bX$ are uniformly bounded as $%
D\to\infty$;

\item[(H3)] $\ds\liminf_{D\to\infty}D^{-1}\lambda_{\min}(\bX_s^{\prime }\bX%
_s)>0$;

\item[(H4)] $\ds\liminf_{D\to\infty}D^{-1}\lambda_{\min}(\cF)>0$.
\end{enumerate}

Theorem \ref{EwEw2} gives an approximation for the second term on the
right-hand side of (\ref{mcpeeb2}). This result uses the additional notation
\begin{align}
&\bx_{dij}=\bx_{di}+\bx_{dj}, \quad \bm_{d} =(\cero_{d-1}^{\prime },1,\cero_{D-d}^{\prime })^{\prime },\quad %
\eeta_{d}=\sigma _{u}^{2}\bV_{s}^{-1}\bZ_{s}\bm_{d},\quad \notag \\
&E_{dij} =\exp \left\{ 2\alpha _{d}+\bx_{dij}^{\prime }\bbeta+\frac{1}{2}\bx%
_{dij}^{\prime }\bQ_{s}\bx_{dij}+2\gamma _{d}\left( \sigma _{u}^{2}-\gamma
_{d}\bar{\bx}_{ds}^{\prime }\bQ_{s}\bar{\bx}_{ds}\right) \right\},
\notag \\
&K_d =\tr\left( \cF^{-1}\frac{\partial \eeta_{d}^{\prime }}{%
\partial \btheta}\bV_{s}\frac{\partial \eeta_{d}}{\partial \btheta}\right)
+\left( \frac{\partial \alpha _{d}}{\partial \btheta}+2\frac{\partial \eeta%
_{d}^{\prime }}{\partial \btheta}\bV_{s}\eeta_{d}\right) ^{\prime }\cF%
^{-1}\left( \frac{\partial \alpha _{d}}{\partial \btheta}+2\frac{\partial %
\eeta_{d}^{\prime }}{\partial \btheta}\bV_{s}\eeta_{d}\right)\notag\\
&M_{2d,ij}(\bbeta,\btheta) =E_{dij} K_d.\notag
\end{align}

\begin{theo}\label{EwEw2}
Let $\hat w_{di}^E=\hat w_{di}(\hat\btheta)$ be the second-stage EB predictor of $w_{di}$, with $\hat\btheta$ denoting either ML or REML estimator of $\btheta$ under the nested-error model with log-transformation (\ref{neemodel}). If assumptions (H1)-(H4) hold, then
\begin{equation*}
E\left\{ (\hat{w}_{di}^{E}-\hat{w}_{di})(\hat{w}_{dj}^{E}-\hat{w}%
_{dj})\right\} =M_{2d,ij}(\bbeta,\btheta)+o(D^{-1}).
\end{equation*}%
\end{theo}


Theorem \ref{CP} gives a second-order unbiased approximation for the first of the crossed product terms in (\ref{mcpeeb2}); the last term is analogous. For this theorem, we need to introduce additional notation. We define
\begin{align}
&E_{dij}^{\ast } =\exp \left\{ \alpha _{d}+\bx_{dij}^{\prime }\bbeta+\sigma
_{e}^{2}+\sigma _{u}^{2}(3+\gamma _{d})+h_{d,ii}+2h_{d,ij}-2\gamma_dh_{d,j}-\gamma_d^2h_d\right\}.  \label{Edij*}
\end{align}%
We also define
$\bE_{d}=2(\bDelta_{1}\eeta_{d},\bDelta_{2}\eeta_{d})$, $\bA%
_{d}=(\alpha _{d,ht})$, with $\alpha _{d,ht}=\partial ^{2}\alpha
_{d}/\partial \theta _{h}\partial \theta _{t}$, $\bB_{d}=(b_{d,ht})$ with $%
b_{d,ht}=2\eeta_{d}^{\prime }\bV_{s}(\partial^2 \eeta_{d}/\partial \theta
_{h}\partial \theta _{t})$,
\begin{equation*}
\bG_{d}=\underset{1\leq k\leq 2}{\col}\left\{ \left( \displaystyle{\frac{%
\partial \alpha _{d}}{\partial \btheta}}+2\displaystyle{\frac{\partial \eeta%
_{d}^{\prime }}{\partial \btheta}}\bV_{s}\eeta_{d}\right) ^{\prime }\cF%
^{-1}\Phi _{k}\right\} ,
\end{equation*}%
for $\Phi _{k}=(\phi _{hk\ell })_{h,\ell }$ with $\phi _{hk\ell }=\tr(\bV_{s}^{-1}\bDelta_{h}\bV_{s}^{-1}\bDelta_{t}%
\bV_{s}^{-1}\bDelta_{k})$, $\bepsilon_{d}=\underset{1\leq h\leq 2}{\col}(4\eeta_{d}^{\prime }\bDelta_{h}%
\eeta_{d})$, $\bvarsigma=(\varsigma_1,\varsigma_2)^{\prime}$, with $\varsigma_h=2\tr(\cF^{-1}\Phi _{h})$, $h=1,2$, and $\bnu=(\nu_1,\nu_2)^{\prime}$, with $\nu_h=\tr(\bP_s\bDelta_h)-\tr(\bV_s^{-1}\bDelta_h)$, $h=1,2$, and
$$
C_d=\tr\left[ \cF^{-1}\left( \frac{\partial \eeta%
_{d}^{\prime }}{\partial \btheta}\bE_{d}+\frac{\bA_{d}+\bB_{d}}{2}-\bG%
_{d}\right) \right] +\left( \frac{\partial \alpha _{d}}{\partial \btheta}+2%
\frac{\partial \eeta_{d}^{\prime }}{\partial \btheta}\bV_{s}\eeta_{d}\right)
^{\prime }\cF^{-1}\left( \bnu+\frac{\bepsilon_{d}+\bvarsigma}{2}\right).
$$
Finally, we define
\begin{align}
& M_{2d,ij}^{\ast }(\bbeta,\btheta)=E_{dij}^{\ast } K_d,\quad T_{d,ij}(\bbeta,\btheta)=E_{dij}C_d,\quad T_{d,ij}^*(\bbeta,\btheta)=E_{dij}^*C_d,\notag\\
& M_{3d,ij}(\bbeta,\btheta)= \frac{1}{2}M_{2d,ij}(\bbeta,\btheta)+T_{d,ij}(\bbeta,\btheta) -\frac{1}{2}M_{2d,ij}^*(\bbeta,\btheta)-T_{d,ij}^*(\bbeta,\btheta).\label{M3dij}
\end{align}


\begin{theo}
\label{CP}  Let $\hat w_{di}^E=\hat w_{di}(\hat\btheta)$ be the second-stage EB predictor of $w_{di}$ under the nested-error model with log-transformation (\ref{neemodel}), with $\hat\btheta$ denoting either ML or REML estimator of $\btheta$. If assumptions (H1)-(H4) hold, then for $i,j\in \bar s_d$, we have
\begin{equation*}
E\left\{(\hat w_{di}^E-\hat w_{di})(\hat w_{dj}-w_{dj})\right\}=M_{3d,ij}(%
\bbeta,\btheta)+o(D^{-1}).
\end{equation*}
If $\hat\btheta$ is the REML estimator, set $\bnu=\cero_2$ in $M_{3d,ij}(\bbeta,\btheta)$.
\end{theo}


Finally, Theorem \ref{FinalMCPE} gives a second-order approximation to the
MCPE of $\hat w_{di}^E$ and $\hat w_{dj}^E$, as a direct consequence of
decomposition (\ref{mcpeeb2}) and Theorems \ref{mseeb1}, \ref{EwEw2} and \ref{CP}.

\begin{theo}
\label{FinalMCPE}  Let $\hat w_{di}^E=\hat w_{di}(\hat\btheta)$ be the second-stage EB predictor of $w_{di}$ under the nested-error model with log-transformation (\ref{neemodel}), with $\hat\btheta$ denoting either ML or REML estimator of $\btheta$. Under assumptions (H1)-(H4), it holds
\begin{equation*}
\mbox{MCPE}(\hat w_{di}^E,\hat w_{dj}^E)=M_{1d,ij}(\bbeta,\btheta)+M_{2d,ij}(%
\bbeta,\btheta)+M_{3d,ij}(\bbeta,\btheta)+M_{3d,ji}(\bbeta,\btheta%
)+o(D^{-1}).
\end{equation*}
\end{theo}

The following corollary gives a second-order approximation to the MSE of the second-stage EB predictor $\hat\tau_d^E$ of the area mean $\tau_d$.
 \begin{coro}\label{FinalMSEtau}
 An approximation to the MSE of $\hat\tau_d^E$ is obtained noting that
 \begin{equation}\label{MSEtau}
 \mbox{MSE}(\hat\tau_d^E)=\frac{1}{N_d^2}\left\{2\sum_{i\in \bar
 s_d}\sum_{j\in \bar s_d,j>i}\mbox{MCPE}(\hat w_{di}^E,\hat w_{dj}^E)+\sum_{i\in \bar
 s_d}\mbox{MSE}(\hat w_{di}^E)\right\}
 \end{equation}
 and applying Theorem \ref{FinalMCPE} to obtain second-order approximations of $\mbox{MCPE}(\hat w_{di}^E,\hat w_{dj}^E)$ and of $\mbox{MSE}(\hat w_{di}^E)$ by setting $i=j$. In fact, going through all the proofs, it can be seen that the remainder term in Theorem \ref{FinalMCPE} is $o(D^{-1})$ uniformly for all $i$ and $j$; in other words,
 $$
 \mbox{MCPE}(\hat w_{di}^E,\hat w_{dj}^E)=\sum_{k=1}^3M_{kd,ij}(\bbeta,\btheta)+M_{3d,ji}(\bbeta,\btheta)+m_{d,ij}(\bbeta,\btheta).
 $$
 where $\max_{1\leq i,j\leq N_d} m_{d,ij}(\bbeta,\btheta)=o(D^{-1})$. This implies that the resulting approximation to the MSE of $\hat\tau_d^E$ is also $o(D^{-1})$.
\end{coro}


\section{Estimation of the uncertainty}

\label{secEstimUncer}

The following theorem states that replacing the unknown parameters $\btheta$
and $\bbeta$ by their corresponding ML estimators $\hat\btheta$ and $\hat%
\bbeta=\tilde\bbeta(\hat\btheta)$ in $M_{1d,ij}(\bbeta,\btheta)$ leads to a $%
O(D^{-1})$ bias. It also gives a second-order approximation for that bias,
which can then be corrected. The proof follows closely that of Theorem 4 in
Molina (2009).

\begin{theo}
\label{them6}  Let $\hat\btheta$ denote either ML or REML estimator of $\btheta$ under the nested-error model with log-transformation (\ref{neemodel}) and $\hat\bbeta=\tilde\bbeta(\hat\btheta)$. If assumptions (H1)-(H4) hold, then
\begin{equation*}
E\{M_{1d,ij}(\hat{\bbeta},\hat{\btheta})\}=M_{1d,ij}(\bbeta,\btheta%
)+\sum_{k=1}^{3}\Lambda _{d,ij,k}(\bbeta,\btheta)+o(D^{-1}),
\end{equation*}%
where
\begin{align*}
\Lambda _{d,ij,1}(\bbeta,\btheta)& =2\left( \partial M_{1d,ij}/\partial %
\btheta\right) ^{\prime }\cF^{-1}\bnu, \\
\Lambda _{d,ij,2}(\bbeta,\btheta)& =(1/2)\tr\left[ \left( \partial
^{2}M_{1d,ij}/\partial \btheta^{2}\right) \cF^{-1}\right] , \\
\Lambda _{d,ij,3}(\bbeta,\btheta)& =M_{1d,ij}(\bbeta,\btheta)\,\bx%
_{dij}^{\prime }(\bX_{s}^{\prime }\bV_{s}^{-1}\bX_{s})^{-1}\bx_{dij}.
\end{align*}
If $\hat\btheta$ is the REML estimator, $\Lambda_{d,ij,1}(\bbeta,\btheta)=0$ because $\bnu=\cero_2$.
\end{theo}

It is not difficult to see that plugging the ML estimators $\hat\btheta$ and
$\hat\bbeta$ for the true values $\btheta$ and $\bbeta$ in the above bias
correction terms leads to negligible bias in the sense
\begin{equation}  \label{Elambda}
E\{\Lambda_{d,ij,k}(\hat\bbeta,\hat\btheta)\}=\Lambda_{d,ij,k}(\bbeta,\btheta%
)+o(D^{-1}),\quad k=1,2,3.
\end{equation}
The same occurs for REML estimators of $\btheta$ and $\bbeta$.
According to Theorem \ref{them6} and equation (\ref{Elambda}), an unbiased estimator
of $\mbox{MCPE}(\tilde w_{di},\tilde w_{dj})$ up to $o(D^{-1})$ terms is
given by
\begin{equation}  \label{mcpebest}
\mbox{mcpe}(\tilde w_{di},\tilde w_{dj})=M_{1d,ij}(\hat\bbeta,\hat\btheta%
)-\sum_{k=1}^3 \Lambda_{d,ij,k}(\hat\bbeta,\hat\btheta).
\end{equation}
Moreover, by Molina (2009), it holds that
\begin{equation}  \label{EM2}
E\{M_{2d,ij}(\hat\bbeta,\hat\btheta)\}=M_{2d,ij}(\bbeta,\btheta)+o(D^{-1}).
\end{equation}
So far we have obtained unbiased estimators up to $o(D^{-1})$ terms of the
first two terms on the right-hand side of (\ref{mcpeeb2}). Thus, in order to
have an unbiased estimator of (\ref{mcpeeb2}) of the same order, it only
remains to estimate unbiasedly $M_{3d,ij}(\bbeta,\btheta)$. The next theorem
states that plugging the ML estimators $\hat\btheta$ and $\hat\bbeta
$ in $M_{3d,ij}(\bbeta,\btheta)$ yields an unbiased estimator of the desired
order.

\begin{theo}
\label{them7} Let $\hat\btheta$ denote either ML or REML estimator of $\btheta$ under the nested-error model with log-transformation (\ref{neemodel}) and $\hat\bbeta=\tilde\bbeta(\hat\btheta)$. If assumptions (H1)-(H4) hold, then
\begin{equation*}
E\{M_{3d,ij}(\hat\bbeta,\hat\btheta)\}=M_{3d,ij}(\bbeta,\btheta)+o(D^{-1}).
\end{equation*}
\end{theo}

The analogous result holds for $M_{3d,ji}(\bbeta,\btheta)=E\{ (\hat{w}%
_{di}-\hat{w}_{di})(\hat{w}_{dj}^{E}-w_{dj})\} +o(D^{-1})$. Finally,
from (\ref{mcpebest}), (\ref{Elambda}) and Theorem \ref{them7}, the
estimator
$$
\mbox{mcpe}(\hat{w}_{di}^{E},\hat{w}_{dj}^{E}) =M_{1d,ij}(\hat{\bbeta},\hat{\btheta})-\sum_{k=1}^{3}\Lambda _{d,ij,k}(\hat{\bbeta},\hat{\btheta})+M_{2d,ij}(\hat{\bbeta},\hat{\btheta})+M_{3d,ij}(\hat{\bbeta},\hat{\btheta})+M_{3d,ji}(\hat{\bbeta},\hat{\btheta})
$$
satisfies
\begin{equation*}
E\{\mbox{mcpe}(\hat{w}_{di}^{E},\hat{w}_{dj}^{E})\}=\mbox{MCPE}(\hat{w}%
_{di}^{E},\hat{w}_{dj}^{E})+o(D^{-1}).
\end{equation*}

\section{Bootstrap estimation of the uncertainty}

\label{SecBoot}

Resampling methods are very popular among practitioners due to their
conceptual simplicity, which also makes them less prone to coding errors.
Under the setup of this paper, the naive bootstrap procedure for finite
populations proposed by Gonz\'alez-Manteiga et al. (2008) can be applied for
the estimation of the MSE of either an individual predictor $\hat w_{di}^E$
or for the predicted area mean $\hat\tau_d^E$. It can also be applied to
estimate the MCPE of two individual predictors $\hat w_{di}^E$ and $\hat
w_{dj}^E$, with $j\neq i$. Here we describe only the steps of the bootstrap
procedure for estimation of the MSE of $\hat\tau_d^E$, because for the other
cases is analogous.

\begin{itemize}
\item[1)] With the available data $(\by_s,\bX_s)$ coming from the sample $s$%
, calculate the ML estimators of the model parameters $\hat\bbeta$ and $\hat%
\btheta=(\hat\sigma_u^2,\hat\sigma_e^2)^{\prime }$.

\item[2)] Generate bootstrap random effects $u_d^*\overset{iid}{\sim }\cN%
(0,\hat\sigma_u^2)$, $d=1,\ldots,D$.

\item[3)] Generate bootstrap errors $e_{di}^*\overset{iid}{\sim }\cN%
(0,\hat\sigma_e^2)$, $i=1,\ldots, N_d$, $d=1,\ldots, D$.

\item[4)] Generate a bootstrap population of response variables from the
fitted model
\begin{equation}  \label{bootmodel}
y_{di}^*=\bx_{di}^{\prime }\hat\bbeta+u_d^*+e_{di}^*,\quad i=1,\ldots, N_d,\
d=1,\ldots, D.
\end{equation}
Let $\tau_d^*=N_d^{-1}\sum_{i=1}^{N_d}\exp(y_{di}^*)$ be the true mean of
area $d$ in this bootstrap population.

\item[5)] Let $\by_{s}^{\ast }$ be the vector with the bootstrap elements whose subscripts are in
the original sample $s$, $\{y_{di}^{\ast };i\in s_{d},d=1,\ldots ,D\}$.
Using the bootstrap sample data $\by_{s}^{\ast }$ and $\bX_{d}$, fit the bootstrap model (\ref{bootmodel}), obtaining new model parameter estimators $\hat{\bbeta}^{\ast }$ and $\hat{\btheta}^{\ast }=(\hat{\sigma}_{u}^{2\ast },\hat{\sigma}%
_{e}^{2\ast })^{\prime }$. Calculate the bootstrap second-stage EB predictor
\begin{equation*}
\hat{\tau}_{d}^{E\ast }=\tilde{\tau}_{d}^{\ast }(\hat{\bbeta}^{\ast },\hat{\btheta}^{\ast })=\frac{1%
}{N_{d}}\left\{ \sum_{i\in s_{d}}\exp (y_{di}^{\ast })+\sum_{i\in \bar{s}%
_{d}}\exp (\tilde{y}_{di}^{\ast }+\hat{\alpha}_{d}^{\ast })\right\} ,
\end{equation*}%
for $\hat{\alpha}_{d}^{\ast }=\{\hat{\sigma}_{u}^{2\ast }(1-\hat{\gamma}%
_{d}^{\ast })\}/2$ and
$\tilde{y}_{di}^{\ast }=\bx_{di}^{\prime }\hat{\bbeta}^{\ast }+\hat{\gamma}%
_{d}^{\ast }(\bar{y}_{ds}^{\ast }-\bar{\bx}_{ds}^{\prime }\hat{\bbeta}^{\ast
})$, where $\bar{y}_{ds}^{\ast }=n_{d}^{-1}\sum_{i\in s_{d}}y_{di}^{\ast }$ and
$\hat{\gamma}_{d}^{\ast }=\hat{\sigma}_{u}^{2\ast }/(\hat{\sigma}_{u}^{2\ast }+\hat{\sigma}_{e}^{2\ast }/n_{d})$.

\item[6)] The bootstrap MSE of $\hat\tau_d^{E*}$ is then
\begin{equation}  \label{bootMSE}
\mbox{MSE}_*(\hat\tau_d^{E*})=E_*(\hat\tau_d^{E*}-\tau_d^*)^2,
\end{equation}
where $E_*$ indicates expectation with respect to the probability
distribution induced by model (\ref{bootmodel}) given the original sample
data $\{y_{di};i\in s_d, d=1,\ldots,D\}$.
\end{itemize}

In practice, (\ref{bootMSE}) is approximated by Monte Carlo, by
repeating Steps 2)--5) a large number of times $B$, and then averaging
over the $B$ replicates. Let $\tau_d^{*(b)}$ be the true parameter in $b$-th
replicate and $\hat\tau_d^{E*(b)}$ be the corresponding second-stage EB
predictor. The Monte Carlo approximation of (\ref{bootMSE}), used here as an
estimator of $\mbox{MSE}(\hat\tau_d^E)$, is given by
\begin{equation}  \label{bootmse}
\mbox{mse}_*(\hat\tau_d^E)=\frac{1}{B}\sum_{b=1}^B(\hat\tau_d^{E*(b)}-%
\tau_d^{*(b)})^2.
\end{equation}
Gonz\'alez-Manteiga et al. (2008) proved the consistency of the bootstrap MSE
of the second-stage EB predictor of a linear parameter by the technique of
imitation. With the available analytical formula for the MCPE given in
Theorem \ref{FinalMCPE}, here the result is analogous. First, by imitating
the proofs of Theorems \ref{mseeb1}, \ref{EwEw2} and \ref{CP} under the bootstrap population
given the original sample data, the bootstrap MCPE can be approximated as
\begin{equation}  \label{bootNmse}
\mbox{MCPE}_*(\hat w_{di}^{E*},\hat w_{dj}^{E*})= \mbox{MCPE}_{N*}(\hat
w_{di}^{E*},\hat w_{dj}^{E*})+o(D^{-1}),
\end{equation}
where
\begin{equation*}
\mbox{MCPE}_{N*}(\hat w_{di}^{E*},\hat w_{dj}^{E*})=M_{1d,ij}(\hat\bbeta,\hat%
\btheta)+M_{2d,ij}(\hat\bbeta,\hat\btheta)+M_{3d,ij}(\hat\bbeta,\hat\btheta%
)+M_{3,j,i}(\hat\bbeta,\hat\btheta).
\end{equation*}
Since ML estimates are consistent and $\mbox{MCPE}_{N*}(\hat
w_{di}^{E*},\hat w_{dj}^{E*})$ is a continuous function of $(\hat\bbeta,\hat%
\btheta)$, then $\mbox{MCPE}_{N*}(\hat w_{di}^{E*},\hat w_{dj}^{E*})$ is
also consistent for $\mbox{MCPE}_N(\hat w_{di}^E,\hat
w_{dj}^E)=M_{1d,ij}(\bbeta,\btheta)+M_{2d,ij}(\bbeta,\btheta)+M_{3d,ij}(%
\bbeta,\btheta)+M_{3,j,i}(\bbeta,\btheta)$. However, due to the presence of
the $O(D^{-1})$ bias terms listed in Theorem \ref{them6}, $\mbox{MCPE}%
_N(\hat w_{di}^E,\hat w_{dj}^E)$ is only first-order and not second-order
unbiased for the true $\mbox{MCPE}(\hat w_{di}^{E},\hat w_{dj}^{E})$, that
is,
\begin{equation*}
E\{\mbox{MCPE}_{N*}(\hat w_{di}^{E*},\hat w_{dj}^{E*})\}=\mbox{MCPE}(\hat
w_{di}^{E},\hat w_{dj}^{E})+O(D^{-1}).
\end{equation*}
For bias corrections of the naive bootstrap estimator (\ref{bootmse}) to achieve a $o(D^{-1})$ bias in the case of
linear parameters, see e.g. Butar and Lahiri (2003) and Pfeffermann and
Tiller (2005). For a bias correction based on double bootstrap, see Hall and Maiti (2006a). These corrections can be directly extended to
estimate our specific non-linear parameters $w_{di}$ or $\tau_d$.
These bias corrections might yield negative MSE estimates. Hall and Maiti (2006b) proposed a positive bias-corrected MSE estimate through double bootstrap, but the second-order unbiasedness property is lost. Thus, ensuring positive MSE estimate and second-order unbiased is still a challenge.

\section{Simulation experiment}\label{sec:simexp}

We carried out a simulation experiment to compare, in terms of bias and MSE under the simple mean model $y_{di}=\mu+u_{d}+e_{di}$, the following estimators of the area means $\tau_d$: (i) the second-stage EB predictor $\hat\tau_d^E$; (ii) the naive predictor $\hat \tau_d^N=N_d^{-1}(\sum_{i\in s_d}w_{di}+\sum_{i\in \bar s_d}\hat w_{di}^N)$, where $\hat w_{di}^N=\exp(\hat y_{di})$; Molina (2009)'s predictor $\hat \tau_d^M=N_d^{-1}(\sum_{i\in s_d}w_{di}+\sum_{i\in \bar s_d}\hat w_{di}^M)$, for $\hat w_{di}^M=\exp(\hat y_{di}+\hat\alpha_d^M)$, with $\hat\alpha_d^M=\hat\sigma_u^2(1-\hat\gamma_d)/2$; (iii) direct estimator $\hat \tau_d^D=n_d^{-1}\sum_{i\in s_d} w_{di}$ and (iv) the estimator obtained assuming the area-level model of Fay and Herriot (1979), $\hat\tau_d^D=\mu+v_d+\epsilon_d$, where  $v_d$ are assumed iid with $E(v_d)=0$, $var(v_d)=\sigma_v^2$, and $\epsilon_d$ are independent with $E(\epsilon_d)=0$ and $var(\epsilon_d)=\psi_d$, with $\psi_d$ assumed to be known and fixed to the sampling variance of the direct estimator $\hat \tau_d^D$, $d=1,\ldots,D$. We will also analyze the contribution of each term of $MSE(\hat\tau_d)$.

\begin{figure}[h]
\begin{minipage}{0.5\linewidth}{
 \centering
 \makebox{\includegraphics[width=75mm,height=75mm]{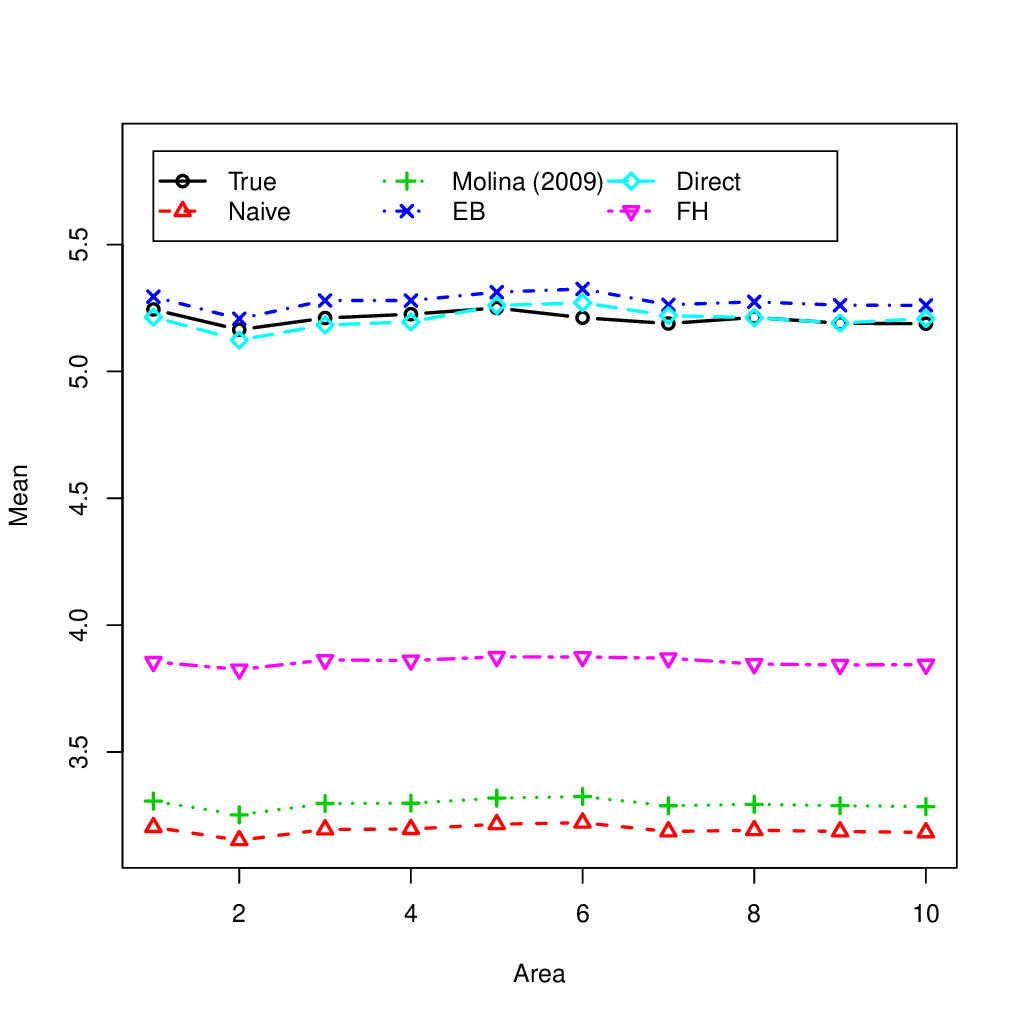}}}\end{minipage}%
\begin{minipage}{0.5\linewidth}{
 \centering
 \makebox{\includegraphics[width=75mm,height=75mm]{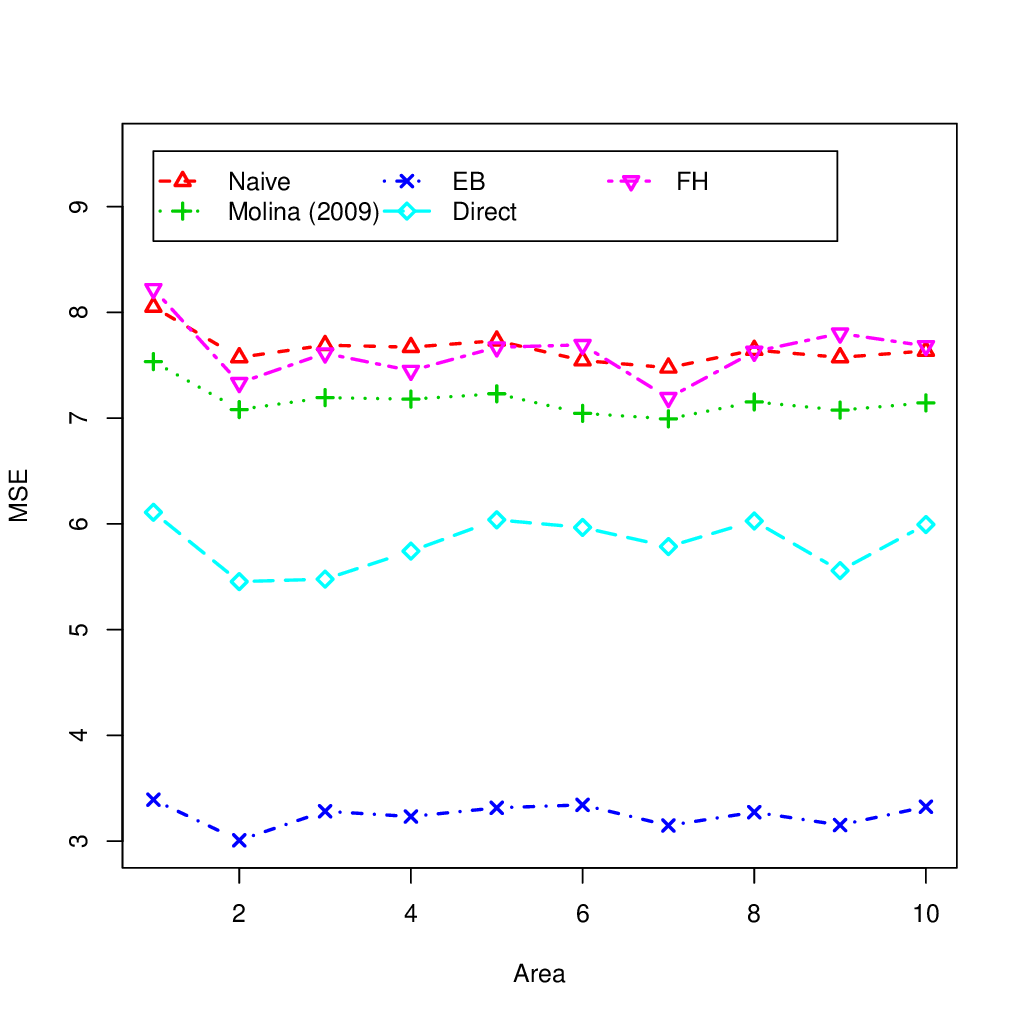}}}\end{minipage}
\vspace{-0.3 cm} 
\caption{Monte Carlo means of true values, naive, Molina (2009) and second-stage EB predictors, direct estimator and estimator based on FH model (left). Monte Carlo MSEs of all the estimators (right).}
\label{BiasAll}
\end{figure}

We consider a limited number of areas, $D=10$, in order to analyze the small sample properties of the estimators. Population sizes of the areas are taken as $N_d=200$, $d=1,\ldots,D$, which gives a total population size of $N=2000$. Model parameters are taken as $\mu=1$, $\sigma_e^2=1$ and $\sigma_u^2=0.3$, leading to a variance fraction $\gamma_d=\sigma_u^2/(\sigma_u^2+\sigma_e^2/n_d)=0.75$. A total of $K=10,000$ Monte Carlo (MC) populations were generated from the mentioned mean model. In each MC simulation replicate, simple random samples $s_d$ without replacement of size $n_d=10$ were drawn independently from each area $d$, making a total sample size of $n=100$. In this case, by Proposition \ref{RBnaiveM}, the actual relative bias of the naive predictor $\tilde w_{di}^N$ amounts to $RB(\tilde w_{di}^N)=-41.6\%$. For Molina (2009)'s predictor, it is $RB(\tilde w_{di}^M)=-39.3\%$. Let us now look at the actual biases and MSEs of each type of estimator of $\tau_d$, $d=1,\ldots,D$. Figure \ref{BiasAll} (left) plots the MC means of the true values $\tau_d$ and of the estimators (i)--(iv) and the MSEs (right). This figure illustrates how the naive and Molina (2009)´s predictors are both considerably biased low and also how the EB predictor $\hat \tau_d^E$ proposed in this paper has a negligible bias together with a substantially smaller MSE than all other estimators.

\begin{figure}[h]
 \centering
 \makebox{\includegraphics[width=100mm,height=75mm]{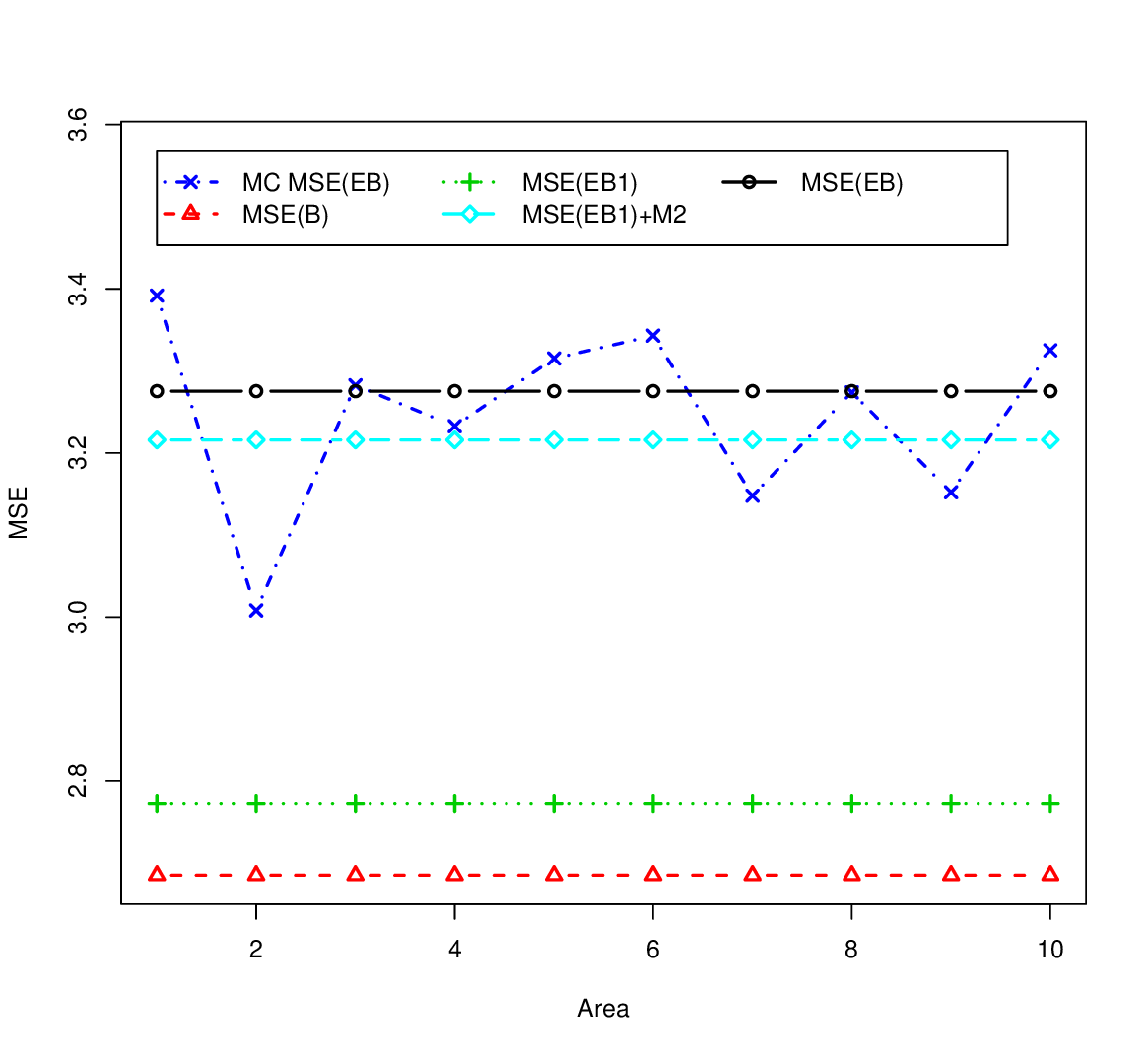}}
\vspace{-0.3 cm} 
\caption{MC MSE of second-stage EB predictor $\hat\tau_d^E$ labelled ``MC MSE(EB)", MSE of best predictor $\tilde\tau_d$ labelled ``MSE(B)", MSE of first-stage EB predictor $\hat\tau_d$ labelled ``MSE(EB1)", the same but adding the crossed-product terms $M_{2d,ij}$, and total MSE of second-stage EB predictor $\hat\tau_d^E$ labelled ``MSE(EB)".}
\label{MSETerms}
\end{figure}

Next we analyze the contribution of each MSE term to the total $MSE(\hat\tau_d^E)$ in this simulation experiment. Figure \ref{MSETerms} displays the
MC approximation to $MSE(\hat\tau_d^E)$ labelled ``MC MSE(EB)", $MSE(\tilde\tau_d)$ given in Corollary \ref{mseebtau} labelled ``MSE(B)", $MSE(\hat\tau_d)$ given in Theorem \ref{mseeb1} labelled ``MSE(EB1)", the same but adding the crossed-product terms $M_{2d,ij}$ given in Theorem \ref{EwEw2}, and finally the analytical approximation to $MSE(\hat\tau_d^E)$ obtained from Theorem \ref{FinalMCPE} and Corollary \ref{FinalMSEtau} that includes the terms $M_{3d,ij}+M_{3d,ji}$. We can clearly see that in this simulation experiment, the naive MSE estimators $MSE(\tilde\tau_d)$ or $MSE(\hat\tau_d)$ underestimate the true MSE to a great extent, and the additional MSE terms of Theorems \ref{EwEw2} and \ref{FinalMCPE} seem to be necessary to avoid undesired underestimation of the MSE.

\section{Estimation of mean income in municipalities from Mexico}

\label{secExample}

In this section we apply the obtained results to the estimation of mean income in municipalities of the State of Mexico. Data comes from two different sources. One is the Module of Socio-economic Conditions (MCS in Spanish) from the 2010 Mexican National Survey on Income and Expense of Households (ENIGH in Spanish). The MCS collects microdata on income, health, nutrition, education, social security, quality of household, basic equipment and social cohesion in Mexico. We also have available micro data from the Census of the same year. The Census contains several of the variables also contained in the MCS, but the income variable used officially (monthly total per capita income) is collected only in the mentioned survey. Based on both data sources, we estimate mean income in each municipality that appears in the MCS survey data (many of them are not sampled by the MCS), except for one which, after a preliminary study of the considered variables, turned out to be very different from the other municipalities (outlier). This makes a total of $D=57$ municipalities. From these, the minimum sample size is 8 and the maximum is 2037, with a median of 96 and an average of 185.

After a preliminary check of the relationships between income and the available variables in the MCS, we selected as auxiliary variables age, $\mbox{age}^2$, $\mbox{age}^3$, the indicators of gender, indigenous population, activity sectors (including unemployed and inactive), composition of household, quality of dwelling, indicator of receiving social benefits, classification according to the available equipment, years of schooling, indicator of rural/urban area and the interactions between quality of dwelling with rural/urban area and of composition of household with gender. Since income distribution in Mexico is highly skewed, the model was fitted to $\log(\mbox{income}+k)$ where $k=171$ was selected to achieve an approximately symmetric distribution of model residuals. The Supplementary material shows the histograms of income before and after the transformation. It shows also the resulting fitted regression parameters. The fitted variance components are $\sigma_u^2=0.0160$ and $\sigma_e^2=0.3245$, which lead to an average contribution of $\sigma_u^2$ to the total variance of $D^{-1}\sum_{d=1}^D\gamma_d=0.789$.

We computed also direct Horvitz-Thompson estimators of mean income $w_{di}$ together with their sampling variances, obtained as
\begin{equation*}
\hat\tau_d^{DIR}=N_d^{-1}\sum_{i\in s_d}\pi_{di}^{-1}w_{di},\quad \mbox{var}%
(\hat\tau_d^{DIR})=N_d^{-2}\sum_{i\in s_d}\pi_{di}^{-2}(1-\pi_{di})w_{di}^2,
\end{equation*}
where $\pi_{di}$ is the inclusion probability of $i$-th unit in the sample from municipality $d$. The sampling variance is obtained using the following approximation for the second-order inclusion probabilities $\pi_{d,ij}\approx\pi_{di}\pi_{dj}$, $j\neq i$, and noting that $\pi_{d,ii}=\pi_{di}$ for all $i$. Figure \ref{Estimates} shows EB, direct, Molina (2009) and naive estimators of mean income for the $D=57$ municipalities. This figure illustrates that direct estimators are somewhat unstable. According to Proposition \ref{RBnaiveM}, Molina (2009) and naive estimators have an average estimated relative bias of -14.77\% and -14.92\% respectively. In this application, both take very similar values (superposed in the plot) and their values are lower than those of EB estimators, which could be due to the mentioned theoretical bias.

\begin{figure}[h]
 \begin{minipage}{0.5\linewidth}{
 \centering
 \makebox{\includegraphics[width=75mm,height=70mm]{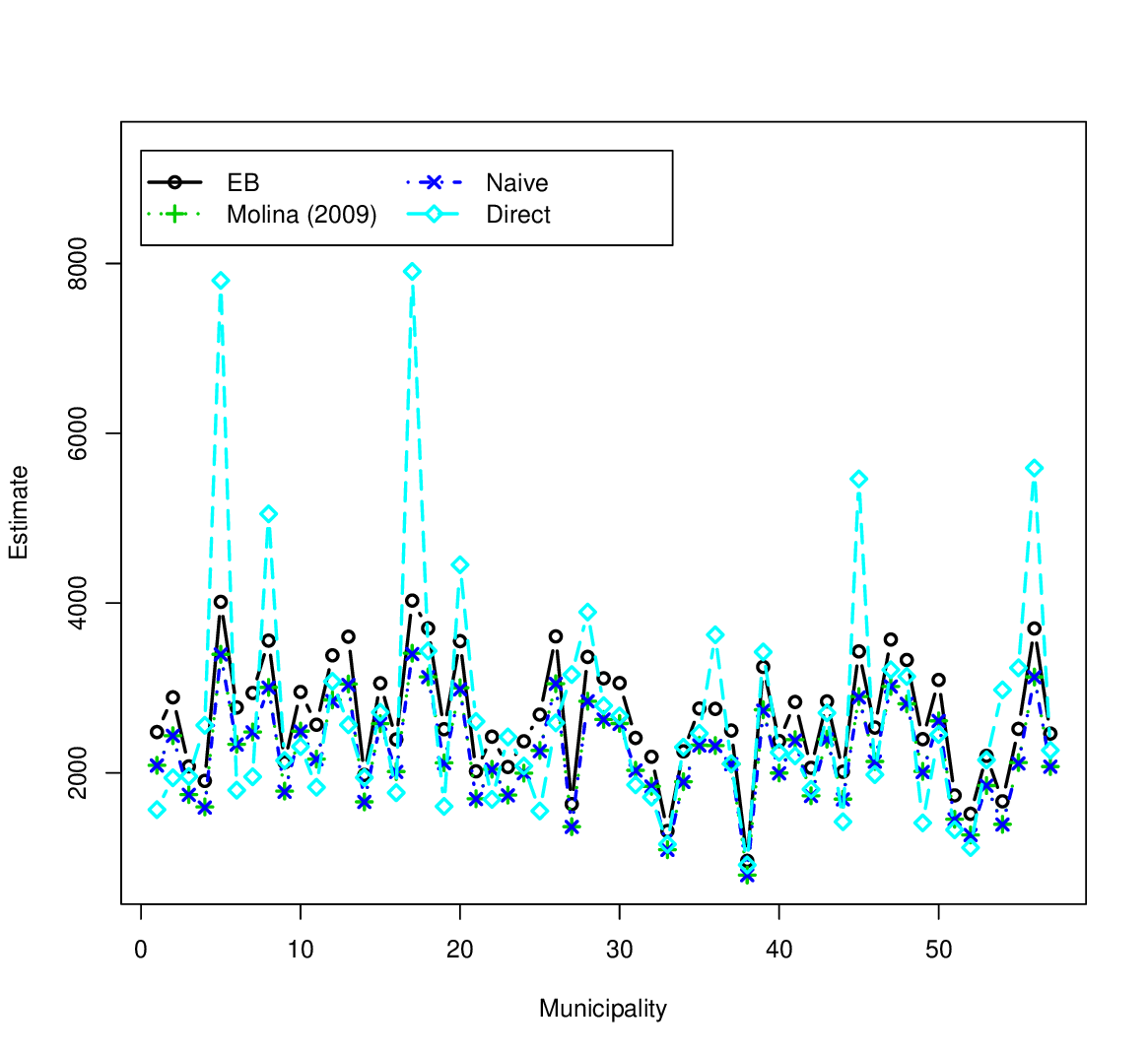}}}\end{minipage}\begin{minipage}{0.5\linewidth}{
 \centering
 \makebox{\includegraphics[width=75mm,height=70mm]{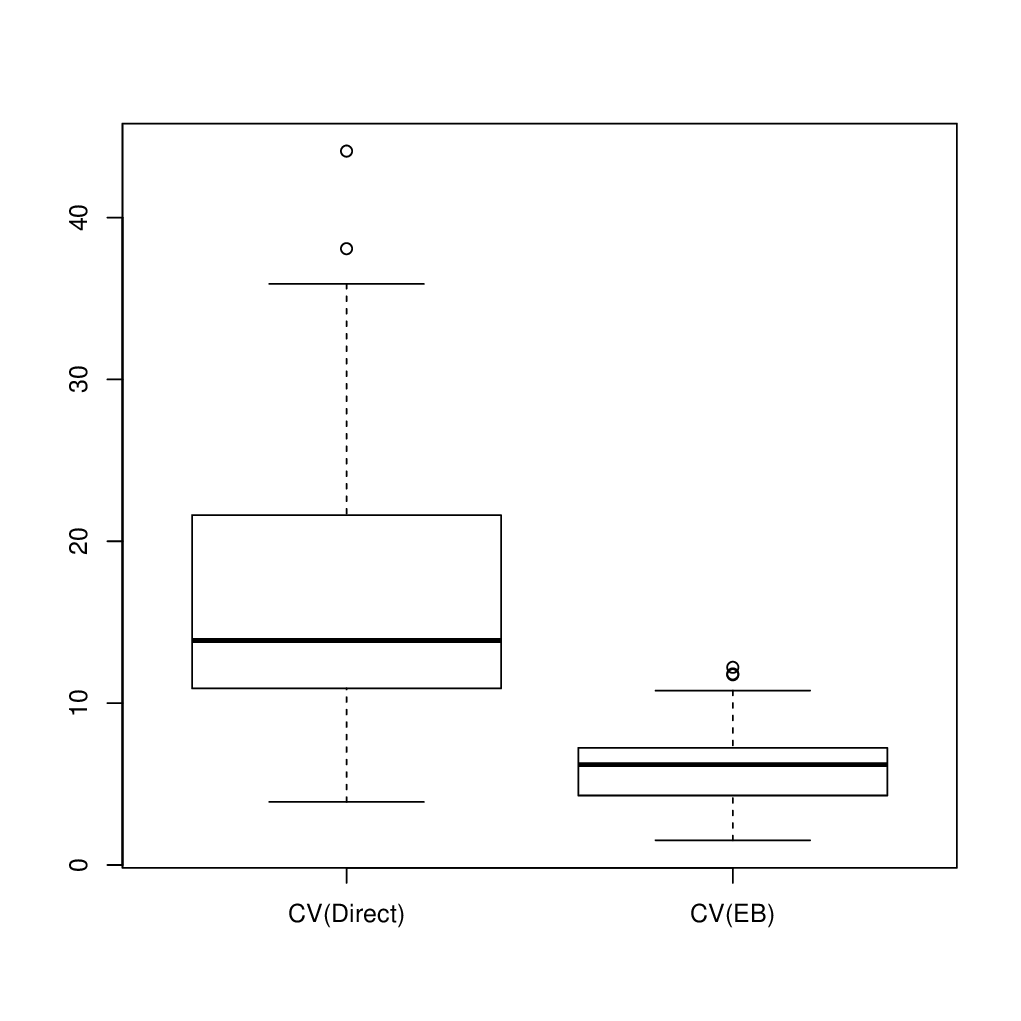}}}\end{minipage}
 \vspace{-0.3 cm}
 \caption{\label{CVs} EB, Molina (2009), naive and direct estimates of mean income for each municipality (left) and boxplots of estimated CVs of direct and EB estimates of mean income (right). }
 \vspace{0.3 cm}
 \end{figure}

Finally, boxplots of the estimated coefficients of variation (CVs) defined for any estimator $\hat\tau_d$ as $\mbox{cv}(\hat\tau_d)=100\times\sqrt{\mbox{var}(\hat\tau_d)}
/\hat\tau_d$ are shown in Figure \ref{CVs}. These boxplots show the significant reduction in CV obtained when using EB estimators instead of the default direct estimators.

\section*{APPENDIX: PROOFS}

In this appendix, the Euclidean norm of a vector $\ba$ is denoted by $|\ba|=(%
\ba^{\prime }\ba)^{1/2}$. For a matrix $A$, we consider the norms $%
\|A\|=\lambda_{\max}^{1/2}(A^{\prime }A)$ and $\|A\|_2=\tr^{1/2}(A^{\prime
}A)$, where $\lambda_{\max}(A)$ denotes the maximum eigenvalue of $A$.
Asymptotic orders refer to $D\to\infty$.

\subsection*{PROOF OF PROPOSITION \protect\ref{RBnaiveM}}

The naive predictor of $w_{di}$ can be expressed in terms of the best predictor $\tilde w_{di}$
as $\tilde w_{di}^N=\exp(\tilde y_{di})=\exp(-\alpha_d)\tilde w_{di}$. Now since the best predictor $\tilde w_{di}$ is unbiased, that is,
$E(\tilde w_{di})=E(w_{di})$, we have $E(\tilde w_{di}^N)=\exp(-\alpha_d)E(\tilde w_{di})=\exp(-\alpha_d)E(w_{di})$. The relative bias (RB) of
$\tilde w_{di}^N$ is then $RB(\tilde w_{di}^N)=E(\tilde w_{di}^N)/E(w_{di})-1=\exp(-\alpha_d)-1$.
Similarly, Molina (2009) predictor can be expressed as $\tilde w_{di}^N=\exp(\tilde y_{di}+\alpha_d^M)=\exp(-\sigma_e^2/2)\tilde w_{di}$. Taking expected value, we get $E(\tilde w_{di}^M)=\exp(-\sigma_e^2/2)E(w_{di})$. Thus, $RB(\tilde w_{di}^M)=E(\tilde w_{di}^M)/E(w_{di})-1=\exp(-\sigma_e^2/2)-1$. \hfill $\Box $


The next lemma is required in the proofs of several of the
remaining results.

\begin{lemm}
\label{lemmaV} Let $\bV_s$ be the covariance matrix of $\by_s$, $\cF$ and $\cF_R$ the ML and REML Fisher-information matrices respectively, and $\bQ_s=(\bX_s^{\prime}\bV_s^{-1}\bX_s)^{-1}$.
It holds

\begin{itemize}
\item[(i)] Condition (H1) implies $\|\bV_s\|=O(1)$.

\item[(ii)] $\|\bV_s^{-1}\|=O(1)$.

\item[(iii)] Conditions (H1) and (H3) imply $\|\bQ_s\|=O(D^{-1})$.

\item[(iv)] Condition (H4) implies $\|\cF^{-1}\|=O(D^{-1})$ and $\|\cF_R^{-1}\|=O(D^{-1})$.
\end{itemize}
\end{lemm}

\subsection*{PROOF OF LEMMA \protect\ref{lemmaV}}

(i) Since $\bV_{s}$ is symmetric and block-diagonal with blocks equal to $\bV%
_{ds}$, $d=1,\ldots ,D$, we have
\begin{equation*}
\Vert \bV_{s}\Vert =\lambda _{\max }^{1/2}(\bV_{s}^{2})=\lambda _{\max }(\bV%
_{s})=\max_{1\leq d\leq D}\{\lambda _{\max }(\bV_{ds})\}.
\end{equation*}%
Now since $\bV_{ds}=\sigma _{u}^{2}\uno_{n_{d}}\uno_{n_{d}}^{\prime }+\sigma
_{e}^{2}\bI_{n_{d}}$, we have
\begin{equation*}
\lambda _{\max }(\bV_{d})\leq \sigma _{u}^{2}\lambda _{\max }(\uno_{n_{d}}%
\uno_{n_{d}}^{\prime })+\sigma _{e}^{2}\lambda _{\max }(\bI_{n_{d}})=\sigma
_{u}^{2}n_{d}+\sigma _{e}^{2}.
\end{equation*}%
Then, by assumption (H1), we obtain
\begin{equation*}
\Vert \bV_{s}\Vert =\max_{1\leq d\leq D}\{\lambda _{\max }(\bV_{ds})\}\leq
\sigma _{u}^{2}\max_{1\leq d\leq D}n_{d}+\sigma _{e}^{2}=O(1),
\end{equation*}%
which implies (i). \hfill $\Box $\newline
(ii) Similarly as before, we have
\begin{equation*}
\Vert \bV_{s}^{-1}\Vert =\lambda _{\max }(\bV_{s}^{-1})=\lambda _{\min
}^{-1}(\bV_{s})=\left\{ \min_{1\leq d\leq D}\lambda _{\min }(\bV%
_{ds})\right\} ^{-1}.
\end{equation*}%
But again, using the expression of $\bV_{ds}=\sigma _{u}^{2}\uno_{n_{d}}\uno%
_{n_{d}}^{\prime }+\sigma _{e}^{2}\bI_{n_{d}}$, we have
\begin{equation*}
\lambda _{\min }(\bV_{d})\geq \sigma _{u}^{2}\lambda _{\min }(\uno_{n_{d}}%
\uno_{n_{d}}^{\prime })+\sigma _{e}^{2}\lambda _{\min }(\bI_{n_{d}})=\sigma
_{e}^{2}>0,
\end{equation*}%
which is true for all $d\in \{1,\ldots ,D\}$ and for all $D$. Therefore, $%
\Vert \bV_{s}^{-1}\Vert =O(1)$. \hfill $\Box $\newline
(iii) By the definition of $\bQ_{s}=(\bX_{s}^{\prime }\bV_{s}^{-1}\bX%
_{s})^{-1}$, we obtain
\begin{equation*}
\Vert \bQ_{s}\Vert =\lambda _{\max }(\bQ_{s})=\lambda _{\min }^{-1}(\bX%
_{s}^{\prime }\bV_{s}^{-1}\bX_{s}).
\end{equation*}%
But by the definition of eigenvalue, we have
\begin{align*}
& \lambda _{\min }(\bX_{s}^{\prime }\bV_{s}^{-1}\bX_{s})=\min_{v}\frac{%
v^{\prime }\bX_{s}^{\prime }\bV_{s}^{-1}\bX_{s}v}{v^{\prime }v}%
=\min_{v}\left( \frac{v^{\prime }\bX_{s}^{\prime }\bV_{s}^{-1}\bX_{s}v}{%
v^{\prime }\bX_{s}^{\prime }\bX_{s}v}\frac{v^{\prime }\bX_{s}^{\prime }\bX%
_{s}v}{v^{\prime }v}\right) \\
& \geq \left( \min_{w}\frac{w^{\prime }\bX_{s}^{\prime }\bV_{s}^{-1}\bX_{s}w%
}{w^{\prime }\bX_{s}^{\prime }\bX_{s}w}\right) \left( \min_{v}\frac{%
v^{\prime }\bX_{s}^{\prime }\bX_{s}v}{v^{\prime }v}\right) =\lambda _{\min }(%
\bV_{s}^{-1})\lambda _{\min }(\bX_{s}^{\prime }\bX_{s}) \\
& =\lambda _{\max }^{-1}(\bV_{s})\lambda _{\min }(\bX_{s}^{\prime }\bX_{s}).
\end{align*}%
Using (i) and assumption (H3), we finally get
\begin{equation*}
D\Vert \bQ_{s}\Vert =D\lambda _{\min }^{-1}(\bX_{s}^{\prime }\bV_{s}^{-1}\bX%
_{s})\leq \frac{\lambda _{\max }(\bV_{s})}{D^{-1}\lambda _{\min }(\bX%
_{s}^{\prime }\bX_{s})}=O(1),
\end{equation*}%
which means that $\Vert \bQ_{s}\Vert =O(D^{-1})$.\hfill $\Box $\newline
 (iv) Condition (H4) implies
 $$
 D\|\cF^{-1}\|=D\lambda_{\max}(\cF^{-1})=\frac{1}{D^{-1}\lambda_{\min}(\cF)}=O(1),
 $$
 which is equivalent to $\|\cF^{-1}\|=O(D^{-1})$. Moreover, note that $\cF=-A/2+B$, where $B=(b_{h\ell})_{h,\ell=1,2}$ and $A=(a_{h\ell})_{h,\ell=1,2}$, for $a_{h\ell}=\tr(\bV_s^{-1}\bDelta_h\bV_s^{-1}\bDelta_{\ell})$ and $b_{h\ell}=\tr(\bP_s\bDelta_h\bP_s\bDelta_{\ell})$, whereas $\cF_R=B/2$. Then,
 \begin{eqnarray*}
 \lambda_{\min}(\cF_R)&=&\lambda_{\min}\left(B-B/2\right)=\lambda_{\min}\left(B-B/2+A/2-A/2\right)\\
 &=&\lambda_{\min}\left\{\cF+(A-B)/2\right\}\geq \lambda_{\min}(\cF)+\frac{1}{2}\lambda_{\min}(A-B).
 \end{eqnarray*}
 But the diagonal elements of $D^{-1}(A-B)$ tend to zero. Indeed
 \begin{eqnarray*}
 b_{hh}-a_{hh}&=&\tr(\bP_s\bDelta_h\bP_s\bDelta_{h})-\tr(\bV_s^{-1}\bDelta_h\bV_s^{-1}\bDelta_{h})\\
 &=&\tr(\bP_s\bDelta_h\bW_s\bDelta_{h})+\tr(\bW_s\bDelta_h\bV_s^{-1}\bDelta_{h}),
 \end{eqnarray*}
 for $\bW_s=\bP_s-\bV_s^{-1}=\bV_s^{-1}\bX_s(\bX_s^{\prime}\bV_s^{-1}\bX_s)^{-1}\bX_s^{\prime}\bV_s^{-1}$. Then,
 $$
 |b_{hh}-a_{hh}|\leq |\tr(\bP_s\bDelta_h\bW_s\bDelta_{h})|+|\tr(\bW_s\bDelta_h\bV_s^{-1}\bDelta_{h})|.
 $$
 Now, for the second term on the right-hand side, we have
 \begin{eqnarray*}
 &&|\tr(\bW_s\bDelta_h\bV_s^{-1}\bDelta_{h})|=|\tr(\bV_s^{-1}\bX_s(\bX_s'\bV_s^{-1}\bX_s)^{-1}\bX_s\bV_s^{-1}\bDelta_h\bV_s^{-1}\bDelta_{h})|\\
 && =|\tr\{(\bX_s'\bV_s^{-1}\bX_s)^{-1/2}\bX_s\bV_s^{-1}\bDelta_h\bV_s^{-1}\bDelta_{h}\bV_s^{-1}\bX_s(\bX_s'\bV_s^{-1}\bX_s)^{-1/2}\}|\\
 && \leq p\lambda_{\max}\{(\bX_s'\bV_s^{-1}\bX_s)^{-1/2}\bX_s\bV_s^{-1}\bDelta_h\bV_s^{-1}\bDelta_{h}\bV_s^{-1}\bX_s(\bX_s'\bV_s^{-1}\bX_s)^{-1/2}\}\\
 && = p\|\bV_s^{-1/2}\bDelta_{h}\bV_s^{-1}\bX_s(\bX_s'\bV_s^{-1}\bX_s)^{-1/2}\|^2\\
 && \leq p\lambda_{\min}^{-4}(\bV_s)\|\bDelta_{h}\|^2\|\bV_s^{-1/2}\bX_s(\bX_s'\bV_s^{-1}\bX_s)^{-1/2}\|^2\\
 && =p\lambda_{\min}^{-4}(\bV_s)\|\bDelta_{h}\|^2=O(1).
 \end{eqnarray*}
 Similarly, it is easy to see that $|\tr(\bP_s\bDelta_h\bW_s\bDelta_{h})|=O(1)$.
 Therefore, it holds that $D^{-1}(a_{hh}-b_{hh})\to 0$ as $D\to\infty$ for $h=1,2$, leading to $\lim_{D\to\infty}\lambda_{\min}\{D^{-1}(A-B)\}=\liminf\lambda_{\min}\{D^{-1}(A-B)\}=0$, which in turn implies
 $$
 \liminf D^{-1}\lambda_{\min}(\cF_R)\geq \liminf D^{-1}\lambda_{\min}(\cF)+\liminf\lambda_{\min}\{D^{-1}(A-B)\}>0.
 $$
 Then, similarly as we did for $\cF$ above, we obtain $\|\cF_R^{-1}\|=O(D^{-1})$.
 \hfill $\Box$

\subsection*{PROOF OF PROPOSITION \protect\ref{proplog}}

First of all, note that
\begin{equation}\label{Ewdi}
E(w_{di})=E[\exp(y_{di})]=\exp\{\bx_{di}^{\prime}\bbeta+(\sigma_u^2+\sigma_e^2)/2\}.
\end{equation}
On the other hand, $\hat{w}_{di}=\exp (\hat{y}_{di}+\alpha
_{d})$, where $\hat{y}_{di}$ is given by
\begin{equation}
\hat{y}_{di}=\bb_{di}^{\prime }\by_{s},  \label{bdiy}
\end{equation}%
for the vector
\begin{equation}
\bb_{di}=\bV_{s}^{-1}\bX_{s}\bQ_{s}\bx_{di}+\sigma _{u}^{2}\bP_{s}\bZ_{s}\bm%
_{d},  \label{bP}
\end{equation}%
where $\bm_{d}=(\cero_{d-1}^{\prime },1,\cero_{D-d}^{\prime })^{\prime }$.
Replacing $\by_{s}=\bX_{s}\bbeta+\bv_{s}$, for $\bv_{s}=%
\bZ_{s}\bu+\be_{s}$ in \eqref{bdiy} and noting that $\bb_{di}^{\prime }\bX_{s}\bbeta=\bx%
_{di}^{\prime }\bbeta$ because $\bP_{s}\bX_{s}=\cero_{n}$, we obtain
$\hat{y}_{di}=\bx_{di}^{\prime }\bbeta+\bb_{di}^{\prime }\bv_{s}$.
Hence, the first-stage EB predictor of $w_{di}$ can be expressed as
\begin{equation}\label{what}
\hat{w}_{di}=\exp (\hat{y}_{di}+\alpha _{d}),\quad \hat{y}_{di}=\bx%
_{di}^{\prime }\bbeta+\bb_{di}^{\prime }\bv_{s}.
\end{equation}%
Taking expected value, we get
\begin{equation}\label{Ehatwdi}
E(\hat{w}_{di})=\exp\left(\bx_{di}^{\prime }\bbeta+\alpha_d+\bb_{di}^{\prime }\bV_{s}\bb_{di}/2\right).
\end{equation}%

Using the definition of $\bb_{di}$ in (\ref{bP}) and $\bP_{s}$ in (%
\ref{pllike}) and taking into account that $\bX_{s}^{\prime }\bP_{s}=\cero%
_{p\times n}$ and $\bP_{s}\bV_s\bP_s=\bP_s$, it is easy to see that
\begin{equation}\label{bbVbb}
\bb_{di}^{\prime }\bV_{s}\bb_{di}=\bx_{di}^{\prime }\bQ_{s}\bx_{di}+(\sigma _{u}^{2})^{2}\bm%
_{d}^{\prime }\bZ_{s}^{\prime }\bV_{s}^{-1}\bZ_{s}\bm_{d} -(\sigma _{u}^{2})^{2}\bm_{d}^{\prime }\bZ_{s}^{\prime }\bV%
_{s}^{-1}\bX_{s}\bQ_{s}\bX_{s}^{\prime }\bV_{s}^{-1}\bZ_{s}\bm_{d}.
\end{equation}%
Since $\bV_{s}=\diag_{1\leq d\leq D}(\bV_{ds})$ with $\bV_{ds}=\sigma
_{u}^{2}\uno_{n_{d}}\uno_{n_{d}}^{\prime }+\sigma _{e}^{2}\bI_{n_{d}}$, $\bm%
_{d}=(\cero_{d-1}^{\prime },1,\cero_{D-d}^{\prime })^{\prime }$, $\bZ_{ds}=%
\diag_{1\leq d\leq D}(\uno_{n_{d}})$ and $\bX_{s}=(\bX_{1s}^{\prime },\ldots
,\bX_{Ds}^{\prime })^{\prime }$, we obtain
\begin{equation}
\bm_{d}^{\prime }\bZ_{s}^{\prime }\bV_{s}^{-1}\bZ_{s}\bm_{d}=\gamma
_{d}/\sigma _{u}^{2},\quad \bm_{d}^{\prime }\bZ_{s}^{\prime }\bV_{s}^{-1}%
\bX_{s}=(\gamma_{d}/\sigma _{u}^{2})\bar{\bx}_{ds}^{\prime }.
\label{mZVZm}
\end{equation}%
Replacing (\ref{mZVZm}) in (\ref{bbVbb}), we finally obtain
\begin{equation}\label{abdiVbdia}
\bb_{di}^{\prime }\bV_{s}\bb_{di} =\bx_{di}^{\prime }\bQ%
_{s}\bx_{di}+\gamma_d(\sigma_u^2-\gamma_d\bar\bx_{ds}^{\prime }\bQ_{s}\bar\bx%
_{ds})=\gamma_d\sigma_u^2+h_{d,ii}-\gamma_d^2h_d,
\end{equation}%
for $h_{d,ii}=\bx_{di}^{\prime }\bQ_{s}\bx_{di}$ and $h_d=\bar\bx_d^{\prime }\bQ_{s}\bx_d$.
Replacing \eqref{abdiVbdia} in \eqref{Ehatwdi}, we obtain
$$
E(\hat{w}_{di})=\exp\left\{\bx_{di}^{\prime }\bbeta+(\sigma_u^2+\sigma_e^2)/2+(h_{d,ii}-\gamma_d^2h_d)/2\right\}.
$$
Finally, by Lemma \ref{lemmaV}, under (H1) and (H3), we have $\|\bQ_s\|=O(D^{-1})$, and using (H2), it holds $|h_{d,ii}|=O(D^{-1})$ and $|h_d|=O(D^{-1})$.
The result then follows by \eqref{bdiy}.\hfill $\Box$

\subsection*{PROOF OF THEOREM \protect\ref{expresBP}}

(i) The best predictor of $w_{di}=\exp (y_{di})$ is equal to $\tilde{w}%
_{di}=E_{\by_{dr}}\{\exp (y_{di})|\by_{ds}\}$. Here we calculate the more
general expectation
$E_{\by_{dr}}\left\{ \exp (\by_{dr}^{\prime }\bb_{d})|\by_{ds}\right\}$,
where $\bb_{d}$ is a non-stochastic vector of size $N_{d}-n_{d}$, $%
d=1,\ldots ,D$. Now using the conditional distribution given in (\ref%
{condist}), this expectation is given by
\begin{equation}
E_{\by_{dr}}\left[ \exp \left( \by_{dr}^{\prime }\bb_{d}\right) |\by_{ds}%
\right] =\exp \left( \bmu_{dr|s}^{\prime }\bb_{d}+\frac{1}{2}\,\bb%
_{d}^{\prime }\bV_{dr|s}\bb_{d}\right) ,  \label{expecEB}
\end{equation}%
because the integral involved is equal to 1. Now (i) follows from the
expressions for $\bmu_{dr|s}$ and $\bV_{dr|s}$ given in (\ref{murs}) and (%
\ref{vrs}), and taking $\bb_{d}$ as a vector with 1 in position $i$ and the
rest of elements equal to zero.\newline
(ii) The best predictor of $\tau _{d}$ is given by
\begin{equation}\label{taudsum}
\tilde{\tau}_{d} =\tilde{\tau}_{d}(\bbeta,\btheta)=E_{\by_{dr}}\left( \tau
_{d}|\by_{ds}\right)  =\frac{1}{N_{d}}\left[ \sum_{i\in s_{d}}\exp (y_{di})+\sum_{i\in \bar{s}%
_{d}}E_{\by_{dr}}\left\{ \exp (y_{di})|\by_{ds}\right\} \right] .
\end{equation}%
The result then follows by straightforward application of (i).\hfill $\Box $


\subsection*{PROOF OF THEOREM \protect\ref{mseeb}}

For $i, j\in \bar s_d$, we need to calculate
\begin{equation}
\mbox{MCPE}(\tilde{w}_{di},\tilde{w}_{dj})=E(\tilde{w}_{di}\tilde{w}_{dj})-E(%
\tilde{w}_{di}w_{dj})-E(w_{di}\tilde{w}_{dj})+E(w_{di}w_{dj}).
\label{decompmse}
\end{equation}%
Since $u_{d}$ and $e_{di}$ are independent for all $i$, the last term on the
right hand side of (\ref{decompmse}) for $i\neq j$ is given by
\begin{equation}
E(w_{di}w_{dj})=\exp \left\{ (\bx_{di}+\bx_{dj})^{\prime }\bbeta\right\}
E\left\{ \exp (2u_{d})\right\} E\left\{ \exp (e_{di})\right\} E\left\{ \exp
(e_{dj})\right\} .  \label{Ewiwj}
\end{equation}%
In contrast, for $i=j$ we have
\begin{equation}
E(w_{di}^{2})=\exp ( 2\bx_{di}^{\prime }\bbeta) E\left\{ \exp
(2u_{d})\right\} E\left\{ \exp (2e_{di})\right\} .  \label{Ewi2}
\end{equation}%
Observe that the expectations appearing on the right hand side of (\ref%
{Ewiwj}) and (\ref{Ewi2}) are respectively the moment generating function
(m.g.f.) of the independent random variables $2u_{d}$, $e_{di}$, $e_{dj}$
and $2e_{di}$, evaluated at $t=1$. Since the m.g.f. of a random variable $%
X\sim \cN(\mu ,\sigma ^{2})$ is given by $M_{X}(t)=\exp (\mu t+\sigma
^{2}t^{2}/2)$, using this expression we get
\begin{equation}
E(w_{di}w_{dj})=\exp \left\{ (\bx_{di}+\bx_{dj})^{\prime }\bbeta+2\sigma
_{u}^{2}+\sigma _{e}^{2}(1+1_{\{i=j\}})\right\} .  \label{Ewij}
\end{equation}%
Now we obtain $E(\tilde{w}_{di}w_{dj})=E\left\{ \exp (\tilde{y}_{di}+\alpha
_{d}+y_{dj})\right\} $. But by model (\ref{neemodel}), we know
\begin{align*}
& y_{dj}=\bx_{dj}^{\prime }\bbeta+u_{d}+e_{dj}, \\
& \tilde{y}_{di}=\bx_{di}^{\prime }\bbeta+\gamma _{d}(\bar{y}_{ds}-\bar{\bx}%
_{ds}^{\prime }\bbeta)=\bx_{di}^{\prime }\bbeta+\gamma _{d}(u_{d}+\bar{e}%
_{ds}),
\end{align*}%
Then,
\begin{equation*}
\tilde{y}_{di}+y_{dj}=(\bx_{di}+\bx_{dj})^{\prime }\bbeta+(1+\gamma
_{d})u_{d}+e_{dj}+\gamma _{d}\bar{e}_{ds}.
\end{equation*}%
Noting that $u_{d}$, $e_{dj}$ for $j\in \bar s_d$ and $\bar e_{ds}$ are independent, we have
\begin{align}
& E(\tilde{w}_{di}w_{dj})=E\left\{ \exp (\tilde{y}_{di}+\alpha
_{d}+y_{dj})\right\}  \label{Eww} \\
&  =\exp \{(\bx_{di}+\bx_{dj})^{\prime }\bbeta\}\exp
(\alpha _{d})E\left[ \exp \left\{ (1+\gamma _{d})u_{d}\right\} \right] E%
\left\{ \exp (e_{dj})\right\}E\left[ \exp\left\{\gamma_{d}\bar{e}_{ds}\right\} \right] .\notag
\end{align}%
Using the m.g.f.'s evaluated at $t=1$ of the random variables involved in (\ref{Eww}), using the expression of $\alpha
_{d}=\frac{1}{2}\{\sigma _{u}^{2}(1-\gamma _{d})+\sigma _{e}^{2}\}$ and the
fact that $\gamma _{d}(\sigma _{u}^{2}+\sigma _{e}^{2}/n_{d})=\sigma
_{u}^{2} $, we get
\begin{equation}
E(\tilde{w}_{di}w_{dj})=\exp \{(\bx_{di}+\bx_{dj})^{\prime }\bbeta+2\sigma
_{u}^{2}+\sigma _{e}^{2}-\sigma_u^2(1-\gamma_d)=E(w_{di}\tilde{w}_{dj})\}.  \label{Etildewij}
\end{equation}%
Finally, we calculate $E(\tilde{w}_{di}\tilde{w}%
_{dj})=E\left\{ \exp (\tilde{y}_{di}+\tilde{y}_{dj}+2\alpha _{d})\right\} $.
Again, by model (\ref{neemodel}), it holds
\begin{equation*}
\tilde{y}_{di}+\tilde{y}_{dj}=(\bx_{di}+\bx_{dj})^{\prime }\bbeta+2\gamma
_{d}(\bar{y}_{ds}-\bar{\bx}_{ds}^{\prime }\bbeta)=(\bx_{di}+\bx%
_{dj})^{\prime }\bbeta+2\gamma _{d}(u_{d}+\bar{e}_{ds}).
\end{equation*}%
Now since
\begin{equation*}
2\gamma _{d}(u_{d}+\bar{e}_{ds})\sim \cN\left\{ 0,4\gamma _{d}^{2}\left(
\sigma _{u}^{2}+\frac{\sigma _{e}^{2}}{n_{d}}\right) \right\} \equiv \cN%
(0,4\gamma _{d}\,\sigma _{u}^{2}),
\end{equation*}%
then using again the m.g.f. of $\gamma _{d}(u_{d}+\bar{e}_{ds})$ evaluated
at $t=1$, we get
\begin{equation*}
E\left[ \exp \{2\gamma _{d}(u_{d}+\bar{e}_{ds})\}\right] =\exp (2\gamma
_{d}\,\sigma _{u}^{2}).
\end{equation*}%
Finally, using the expression of $\alpha _{d}=\{\sigma
_{u}^{2}(1-\gamma _{d})+\sigma _{e}^{2}\}/2$, we get
\begin{equation}
E(\tilde{w}_{di}\tilde{w}_{dj})=\exp \{(\bx_{di}+\bx_{dj})^{\prime }\bbeta%
\}\exp \{2\sigma _{u}^{2}+\sigma _{e}^{2}-\sigma _{u}^{2}(1-\gamma _{d})\}=E(\tilde{w}_{di}w_{dj}).
\label{Ettildewij}
\end{equation}%
The result follows by replacing (\ref{Ewij}), (\ref{Etildewij}) and (\ref%
{Ettildewij}) in (\ref{decompmse}).\hfill $\Box $


\subsection*{PROOF OF COROLLARY \protect\ref{mseebtau}}

The MSE of $\tilde\tau_d$ is given by
\begin{equation*}
\mbox{MSE}(\tilde\tau_d)=\frac{1}{N_d^2}\left\{2\sum_{i\in \bar
s_d}\sum_{j\in \bar s_d,j>i}\mbox{MCPE}(\tilde w_{di},\tilde
w_{dj})+\sum_{i\in \bar s_d}\mbox{MSE}(\tilde w_{di})\right\}.
\end{equation*}
The result follows by using Theorem \ref{mseeb} separately for $i\neq j$ and
for $i=j$.\hfill$\Box$


\subsection*{PROOF OF THEOREM \protect\ref{mseeb1}}

The mean crossed product error of a pair of individual first-stage
predictors $\hat{w}_{di}$ and $\hat{w}_{dj}$, for $i,j\in \bar s_d$, is given by
\begin{equation}
MCPE(\hat{w}_{di},\hat{w}_{dj})=E(\hat{w}_{di}\hat{w}%
_{dj})+E(w_{di}w_{dj})-E(\hat{w}_{di}w_{dj})-E(w_{di}\hat{w}_{dj}).
\label{MCPEij}
\end{equation}%
The second term on the right hand side of (\ref{MCPEij}) is given in (\ref%
{Ewij}). Concerning the first term on the right hand side of (\ref{MCPEij}),
see that for all $i\in \bar{s}_{d}$, using \eqref{what}, we get
\begin{equation*}
E(\hat{w}_{di}\hat{w}_{dj})=\exp \{2\alpha _{d}+(\bx_{di}+\bx_{dj})^{\prime }%
\bbeta\}E\left[ \exp \{(\bb_{di}+\bb_{dj})^{\prime }\bv_{s}\}\right] ,
\end{equation*}%
where the expectation on the right hand side is the m.g.f. of the normal
random vector $(\bb_{di}+\bb_{dj})^{\prime }\bv_{s}$ evaluated at 1, that
is,
\begin{equation}
E(\hat{w}_{di}\hat{w}_{dj})=\exp \left\{ (\bx_{di}+\bx_{dj})^{\prime }\bbeta+%
(\bb_{di}+\bb_{dj})^{\prime }\bV_{s}(\bb_{di}+\bb_{dj})/2+2\alpha
_{d}\right\} .  \label{Ehatwij}
\end{equation}

Concerning the remaining expectations in (\ref{MCPEij}), noting that
$w_{di}=\exp (y_{di})$ for $y_{di}=\bx_{di}^{\prime }\bbeta+u_d+e_{di}$ and using
(\ref{what}), we can write
\begin{equation*}
y_{di}+\hat{y}_{dj}=(\bx_{di}+\bx_{dj})^{\prime }\bbeta+\bb_{dj}^{\prime }\bv_{s}+u_d+e_{di}+\alpha_d.
\end{equation*}%
Replacing now $\bv_{s}=\bZ_s\bu+\be_s$ and writing $u_d=\bm_d'\bu$, we obtain
\begin{equation}
E(w_{di}\hat{w}_{dj})=\exp \left\{ (\bx_{di}+\bx_{dj})^{\prime }\bbeta+\alpha_{d}\right\}
E\{\exp(e_{di})\}E\left[\exp\{(\bm_d'+\bb_{dj}'\bZ_s)\bu\}\right]E\{\exp(\bb_{dj}'\be_s)\}.
\end{equation}%
Similarly as before, using the m.g.f. of the normal random
vectors involved in the previous expression and rearranging the terms, we obtain
\begin{equation}
E(w_{di}\hat{w}_{dj}) = \exp \left\{ (\bx_{di}+\bx_{dj})^{\prime }\bbeta+\alpha_d+(\sigma_e^2+\sigma_u^2)/2+\bb_{dj}'\bV_s\bb_{dj}/2+\sigma_u^2\bm_d'\bZ_s\bb_{dj}\right\}.  \label{Ewihatwj}
\end{equation}%
Replacing (\ref{Ewij}), (\ref{Ehatwij}) and (\ref{Ewihatwj}) in (\ref{MCPEij}%
), we get
\begin{align}
& \hspace{-0.7cm}\mbox{MCPE}(\hat{w}_{di},\hat{w}_{dj})=\exp \left\{ (\bx%
_{di}+\bx_{dj})^{\prime }\bbeta\right\} \left[ \exp \left\{ 2\sigma
_{u}^{2}+\sigma _{e}^{2}(1+1_{\{i=j\}})\right\} \right.   \label{MCPE1} \\
& +\exp \left\{ (\bb_{di}+\bb_{dj})^{\prime }\bV_{s}(\bb_{di}+\bb%
_{dj})/2+2\alpha _{d}\right\}   \notag \\
& -\exp \left\{(\sigma_e^2+\sigma_u^2)/2+\bb_{di}'\bV_s\bb_{di}/2+\sigma_u^2\bm_d'\bZ_s\bb_{di}+\alpha_d\right\}  \notag \\
& \left. - \exp \left\{ (\sigma_e^2+\sigma_u^2)/2+\bb_{dj}'\bV_s\bb_{dj}/2+\sigma_u^2\bm_d'\bZ_s\bb_{dj}\right\}+\alpha_d\right] .  \notag
\end{align}

Let us calculate the expression of each term in \eqref{MCPE1}. Now similarly as in \eqref{abdiVbdia}, using the definition of $\bb_{di}$ given in (\ref{bP}) and $\bP_{s}$ in (\ref{pllike}), we get
\begin{equation}
(\bb_{di}+\bb_{dj})^{\prime }\bV_{s}(\bb_{di}+\bb_{dj})=(\bx_{di}+\bx%
_{dj})^{\prime }\bQ_{s}(\bx_{di}+\bx_{dj})+4\gamma _{d}\left( \sigma
_{u}^{2}-\gamma _{d}\bar{\bx}_{ds}^{\prime }\bQ_{s}\bar{\bx}_{ds}\right) .
\label{bbVbb2}
\end{equation}
On the other hand, using \eqref{mZVZm}, we get
\begin{equation}\label{sumZb}
\sigma_u^2\bm_d'\bZ_s\bb_{di}=\gamma_d\left(\sigma_u^2+\bx_{di}'\bQ_s\bar\bx_{ds}-\gamma_d\bar\bx_{ds}^{\prime }\bQ_{s}\bar\bx%
_{ds}\right).
\end{equation}
Replacing \eqref{bbVbb2}, \eqref{abdiVbdia}, \eqref{sumZb} and the expression for $\alpha_d$ in \eqref{MCPE1}, we obtain the desired expression for
$\mbox{MCPE}(\hat{w}_{di},\hat{w}_{dj})$.\hfill $%
\Box $

\subsection*{PROOF OF THEOREM \protect\ref{EwEw2}}

We prove it for the case in which $\hat\btheta$ is the ML estimator of $\btheta$. For the REML estimator the proof is analogous, but in fact simpler. Following the same arguments as in the proof of Theorem 1 in Molina (2009),
we obtain
\begin{equation}
E\left\{ (\hat{w}_{di}^{E}-\hat{w}_{di})(\hat{w}_{dj}^{E}-\hat{w}%
_{dj})\right\} =E\left\{ \left( \bh_{di}^{\prime }\cF^{-1}\bs\right) \left( %
\bh_{dj}^{\prime }\cF^{-1}\bs\right) \right\} +o(D^{-1}),  \label{EwEwE}
\end{equation}%
where $\bh_{di}=\partial \hat{w}_{di}/\partial \btheta$. Using the
same ideas as in Theorem 2 in Molina (2009), we get
\begin{align}
& \hspace{-0.7cm}E\left\{ \left( \bh_{di}^{\prime }\cF^{-1}\bs\right) \left( %
\bh_{dj}^{\prime }\cF^{-1}\bs\right) \right\} = \label{hFs2}\\
& \hspace{-0.7cm}\exp \left\{ 2\alpha _{d}+(\bx_{di}+\bx_{dj})^{\prime }%
\bbeta+\frac{1}{2}(\bb_{di}+\bb_{dj})^{\prime }\bV_{s}(\bb_{di}+\bb%
_{dj})\right\} \left\{ \tr\left( \cF^{-1}\frac{\partial \eeta_{d}^{\prime }}{%
\partial \btheta}\bV_{s}\frac{\partial \eeta_{d}}{\partial \btheta}\right)
\right.  \notag \\
& \hspace{-0.7cm}\left. +\left( \frac{\partial \eeta_{d}^{\prime }}{\partial %
\btheta}\bV_{s}(\bb_{di}+\bb_{dj})+\frac{\partial \alpha _{d}}{\partial %
\btheta}\right) ^{\prime }\cF^{-1}\left( \frac{\partial \eeta_{d}^{\prime }}{%
\partial \btheta}\bV_{s}(\bb_{di}+\bb_{dj})+\frac{\partial \alpha _{d}}{%
\partial \btheta}\right) \right\} +o(D^{-1}).  \notag
\end{align}%
Note that by (\ref{bP}), we can express $\bb_{di}$ in terms of $\eeta_{d}$
as follows
\begin{equation}
\bb_{di}=\eeta_{d}+\bV_{s}^{-1}\bX_{s}\bQ_{s}(\bx_{di}-\bX_{s}^{\prime }\eeta%
_{d}).  \label{bbeta}
\end{equation}%
But $\Vert \bZ_{s}\Vert =O(1)$ by assumption (H1). Moreover, $|\bm%
_{d}|=1$. Using Lemma \ref{lemmaV} (ii), we get
\begin{equation}
|\eeta_{d}|=\sigma _{u}^{2}|\bV_{s}^{-1}\bZ_{s}\bm_{d}|\leq \sigma
_{u}^{2}\Vert \bV_{s}^{-1}\Vert \Vert \bZ_{s}\Vert |\bm_{d}|=O(1).
\label{normeta}
\end{equation}%
Now observe that by Lemma \ref{lemmaV} (iii), we have
\begin{equation*}
\Vert \bV_{s}^{-1/2}\bX_{s}\bQ_{s}\Vert =\lambda _{\max }^{1/2}(\bQ_{s}\bX%
_{s}\bV_{s}^{-1}\bX_{s}\bQ_{s})=\lambda _{\max }^{1/2}(\bQ_{s})=O(D^{-1/2}).
\end{equation*}%
Since $\bX_{s}^{\prime }\eeta_{d}=\bX_{ds}^{\prime }\bV_{ds}^{-1}\uno%
_{n_{d}}$, which has bounded norm, and $|\bx_{di}-\bX_{s}^{\prime }\eeta%
_{d}|\leq |\bx_{di}|+|\bX_{s}^{\prime }\eeta_{d}|$, by assumptions
(H1)-(H3), we have
\begin{equation}
|\bV_{s}^{-1}\bX_{s}\bQ_{s}(\bx_{di}-\bX_{s}^{\prime }\eeta_{d})|\leq \Vert %
\bV_{s}^{-1/2}\Vert \Vert \bV_{s}^{-1/2}\bX_{s}\bQ_{s}\Vert |\bx_{di}-\bX%
_{s}^{\prime }\eeta_{d}|=O(D^{-1/2}).  \label{VXQx}
\end{equation}%
From (\ref{bbeta}), (\ref{normeta}) and (\ref{VXQx}), we have obtained
\begin{equation}
\bb_{di}=\eeta_{d}+\ff_{di},\quad |\eeta_{d}|=O(1),\quad |\ff%
_{di}|=O(D^{-1/2}).  \label{decompbdi}
\end{equation}%
Note also that $|\partial \eeta_{d}/\partial \theta _{h}|=O(1)$, since
\begin{equation*}
\frac{\partial \eeta_{d}}{\partial \theta _{h}}=\bV_{s}^{-1}\left( \frac{%
\partial \sigma _{u}^{2}}{\partial \theta _{h}}\bI_{n}-\bDelta_{h}\bV%
_{s}^{-1}\right) \bZ_{s}\bm_{d},\quad h=1,2.
\end{equation*}%
This implies $\left\Vert \partial \eeta_{d}/\partial \btheta\right\Vert
=O(1) $, because
\begin{equation*}
\left\Vert \frac{\partial \eeta_{d}}{\partial \btheta}\right\Vert \leq
\left\Vert \frac{\partial \eeta_{d}}{\partial \btheta}\right\Vert _{2}=\tr%
^{1/2}\left\{ \left( \frac{\partial \eeta_{d}}{\partial \btheta}\right)
^{\prime }\frac{\partial \eeta_{d}}{\partial \btheta}\right\} =\left(
\sum_{h=1}^{2}\left\vert \frac{\partial \eeta_{d}}{\partial \theta _{h}}%
\right\vert ^{2}\right) ^{1/2}\leq 2^{1/2}\max_{h\in \{1,2\}}\left\vert
\frac{\partial \eeta_{d}}{\partial \theta _{h}}\right\vert
\end{equation*}%
By (\ref{bbeta}) and (\ref{VXQx}), we get for any $i$,
\begin{equation}
\cF^{-1}\frac{\partial \eeta_{d}^{\prime }}{\partial \btheta}\bV_{s}\bb_{di}=%
\cF^{-1}\frac{\partial \eeta_{d}^{\prime }}{\partial \btheta}\bV_{s}\eeta%
_{d}+\bkappa_{di},\quad |\bkappa_{di}|=o(D^{-1}).  \label{FetaV}
\end{equation}%
Using repeatedly (\ref{FetaV}), we obtain
\begin{align*}
& \left\{ \frac{\partial \eeta_{d}^{\prime }}{\partial \btheta}\bV_{s}(\bb%
_{di}+\bb_{dj})+\frac{\partial \alpha _{d}}{\partial \btheta}\right\}
^{\prime }\cF^{-1}\left\{ \frac{\partial \eeta_{d}^{\prime }}{\partial \btheta%
}\bV_{s}(\bb_{di}+\bb_{dj})+\frac{\partial \alpha _{d}}{\partial \btheta}%
\right\} \\
& =\left( 2\frac{\partial \eeta_{d}^{\prime }}{\partial \btheta}\bV_{s}\eeta%
_{d}+\frac{\partial \alpha _{d}}{\partial \btheta}\right) ^{\prime }\cF%
^{-1}\left( 2\frac{\partial \eeta_{d}^{\prime }}{\partial \btheta}\bV_{s}%
\eeta_{d}+\frac{\partial \alpha _{d}}{\partial \btheta}\right) +o(D^{-1})
\end{align*}%
and using (\ref{bbVbb2}), we obtain
\begin{equation}
\exp \left\{ 2\alpha _{d}+(\bx_{di}+\bx_{dj})^{\prime }\bbeta+(\bb%
_{di}+\bb_{dj})^{\prime }\bV_{s}(\bb_{di}/2+\bb_{dj})\right\} =E_{dij}.
\label{proofEdij}
\end{equation}%
Replacing (\ref{proofEdij}) in (\ref%
{hFs2}) and then (\ref{hFs2}) in (\ref{EwEwE}), we get the desired
result. \hfill $\Box $


\subsection*{PROOF OF THEOREM \protect\ref{CP}}

Again, we show the result for the ML estimator $\hat\btheta$ of $\btheta$, because for REML the proof is analogous but simpler. The proof is based on the following chain of results:
\begin{itemize}
\item[\textbf{(A)}] For every $\nu\in (0,1)$, there exists a subset of the
sample space $\cB$ on which, for large $D$, it holds
\begin{equation*}
\hat w_{di}^E-\hat w_{di}=\bh_{di}^{\prime }\cF^{-1}\bs+\bh_{di}^{\prime }\cF%
^{-1}(H+\cF)\cF^{-1}\bs+\frac{1}{2}\bh_{di}^{\prime }\cF^{-1}\bd+\frac{1%
}{2}\bs^{\prime }\cF^{-1}S_{di}\cF^{-1}\bs+\br_{di},
\end{equation*}
where $\bh_{di}=\partial \hat w_{di}/\partial\btheta$, $S_{di}=\partial^2
\hat w_{di}/\partial\btheta^2$, $\bd=(d_1,d_2)^{\prime }$, with $d_h=\bs%
^{\prime }\cF^{-1}D_h\cF^{-1}\bs$, $D_h=\partial H/\partial\theta_h$, $h=1,2$%
, and the remainder term $\br_{di}$ satisfies $|\br_{di}|<D^{-3\nu/2}w$, for
a random variable $w$ with bounded first and second moments.

\item[\textbf{(B)}] If $1_{\cB}$ is the indicator function of the set $\cB$,
it holds that
\begin{align}
& \hspace{-0.8cm}E\left\{ (\hat{w}_{di}^{E}-\hat{w}_{di})(\hat{w}%
_{dj}-w_{dj})1_{\cB}\right\} =E\left\{ \bh_{di}^{\prime }\cF^{-1}\bs(\hat{w}%
_{dj}-w_{dj})1_{\cB}\right\}   \label{DecompCP} \\
& \hspace{-0.8cm}\qquad +E\left\{ \bh_{di}^{\prime }\cF^{-1}(H+\cF)\cF%
^{-1}\bs(\hat{w}_{dj}-w_{dj})1_{\cB}\right\}  \notag \\
& \hspace{-0.8cm}\qquad +E\left\{ \frac{1}{2}\bh_{di}^{\prime }\cF^{-1}\bd(%
\hat{w}_{dj}-w_{dj})1_{\cB}\right\}  \notag \\
& \hspace{-0.8cm}\qquad +E\left\{ \frac{1}{2}\bs^{\prime }\cF^{-1}S_{di}\cF%
^{-1}\bs(\hat{w}_{dj}-w_{dj})1_{\cB}\right\} +o(D^{-1}). \notag
\end{align}

\item[\textbf{(C)}] $E\left\{ (\hat{w}_{di}^{E}-\hat{w}_{di})(\hat{w}%
_{dj}-w_{dj})1_{\cB^{c}}\right\} =o(D^{-1})$.

\item[\textbf{(D)}] It holds that
\begin{align}
& \hspace{-0.8cm}E\left\{ \bh_{di}^{\prime }\cF^{-1}\bs(\hat{w}%
_{dj}-w_{dj})\right\} +E\left\{ \bh_{di}^{\prime }\cF^{-1}(H+\cF)\cF%
^{-1}\bs(\hat{w}_{dj}-w_{dj})\right\}  \label{termssumM} \\
& \hspace{-0.8cm}\qquad +E\left\{ \frac{1}{2}\bh_{di}^{\prime }\cF^{-1}\bd(%
\hat{w}_{dj}-w_{dj})\right\} +E\left\{ \frac{1}{2}\bs^{\prime }\cF^{-1}S_{di}%
\cF^{-1}\bs(\hat{w}_{dj}-w_{dj})\right\}  \notag \\
& \hspace{-0.8cm}\qquad =M_{3d,ij}(\bbeta,\btheta)+o(D^{-1}).\notag
\end{align}

\item[\textbf{(E)}] It holds that
\begin{align*}
& \hspace{-0.8cm}E\left\{ \bh_{di}^{\prime }\cF^{-1}\bs(\hat{w}%
_{dj}-w_{dj})1_{\cB^{c}}\right\} =o(D^{-1}), \\
& \hspace{-0.8cm}E\left\{ \bh_{di}^{\prime }\cF^{-1}(H+\cF)\cF^{-1}\bs(%
\hat{w}_{dj}-w_{dj})1_{\cB^{c}}\right\} =o(D^{-1}), \\
& \hspace{-0.8cm}E\left\{ \frac{1}{2}\bh_{di}^{\prime }\cF^{-1}\bd(\hat{w}%
_{dj}-w_{dj})1_{\cB^{c}}\right\} =o(D^{-1}), \\
& \hspace{-0.8cm}E\left\{ \frac{1}{2}\bs^{\prime }\cF^{-1}S_{di}\cF^{-1}\bs(%
\hat{w}_{dj}-w_{dj})1_{\cB^{c}}\right\} =o(D^{-1}).
\end{align*}
\end{itemize}

Applying in turn (C) and (B), we obtain
\begin{eqnarray}
&&\hspace{-0.8cm}E\left\{ (\hat{w}_{di}^{E}-\hat{w}_{di})(\hat{w}%
_{dj}-w_{dj})\right\}   \notag \\
&&\hspace{-0.8cm}\quad =E\left\{ \bh_{di}^{\prime }\cF^{-1}\bs(\hat{w}%
_{dj}-w_{dj})1_{\cB}\right\} +E\left\{ \bh_{di}^{\prime }\cF^{-1}(H+\cF)%
\cF^{-1}\bs(\hat{w}_{dj}-w_{dj})1_{\cB}\right\}   \notag \\
&&\hspace{-0.8cm}\quad +E\left\{ \frac{1}{2}\bh_{di}^{\prime }\cF^{-1}\bd(%
\hat{w}_{dj}-w_{dj})1_{\cB}\right\} +E\left\{ \frac{1}{2}\bs^{\prime }\cF%
^{-1}S_{di}\cF^{-1}\bs(\hat{w}_{dj}-w_{dj})1_{\cB}\right\} +o(D^{-1}).
\notag
\end{eqnarray}%
Finally, writing $1_{\cB}=1-1_{\cB^{c}}$ and applying (E) and (D), we obtain
\begin{equation}
\hspace{-0.8cm}E\left\{ (\hat{w}_{di}^{E}-\hat{w}_{di})(\hat{w}%
_{dj}-w_{dj})\right\} =M_{3d,ij}(\bbeta,\btheta)+o(D^{-1}).  \notag
\end{equation}%
Next we give the proofs of results (A)--(E). \vspace{0.3cm}\newline
\textbf{Proof of (A):} It is obtained by applying Lemma 3 of Molina (2009)
to $\hat{w}_{di}^{E}=\hat{w}_{di}(\btheta)$, where $\hat{\btheta}$ is the ML
estimator of $\btheta$. \vspace{0.3cm}\newline
\textbf{Proof of (B):} Applying (A) we obtain
\begin{align*}
& E\left\{ (\hat{w}_{di}^{E}-\hat{w}_{di})(\hat{w}_{dj}-w_{dj})1_{\cB%
}\right\} =E\left\{ \bh_{di}^{\prime }\cF^{-1}\bs(\hat{w}_{dj}-w_{dj})1_{\cB%
}\right\}  \\
& \qquad +E\left\{ \bh_{di}^{\prime }\cF^{-1}(H+\cF)\cF^{-1}\bs(\hat{w}%
_{dj}-w_{dj})1_{\cB}\right\} +E\left\{ \frac{1}{2}\bh_{di}^{\prime }\cF^{-1}%
\bd(\hat{w}_{dj}-w_{dj})1_{\cB}\right\}  \\
& \qquad +E\left\{ \frac{1}{2}\bs^{\prime }\cF^{-1}S_{di}\cF^{-1}\bs(\hat{w}%
_{dj}-w_{dj})1_{\cB}\right\} +E\left\{ \br_{di}(\hat{w}_{dj}-w_{dj})1_{\cB%
}\right\} .
\end{align*}%
But by Theorem \ref{mseeb1}, we know that $\mbox{MSE}(\hat{w}_{dj})=O(1)$ as
$D$ tends to infinity. Then, applying Hölder's inequality and taking $\nu
\in (2/3,1)$, we obtain
\begin{align}
E\left\{ \br_{di}(\hat{w}_{dj}-w_{dj})1_{\cB}\right\} & \leq E^{1/2}(\br%
_{di}^{2}1_{\cB})E^{1/2}\{(\hat{w}_{dj}-w_{dj})^{2}\}  \label{remain} \\
& <D^{-3\nu /2}E^{1/2}(w^{2})\{\mbox{MSE}(\hat{w}_{dj})\}^{1/2}=o(D^{-1}).\notag
\end{align}%
\textbf{Proof of (C):} Noting that $\hat{w}_{di}^{E}=\exp (\hat{y}_{di}^{E}+%
\hat{\alpha}_{d})$, for $\hat{y}_{di}^{E}=\hat{y}_{di}(\hat{\btheta})$ and $%
\hat{\alpha}_{d}=\alpha _{d}(\hat{\btheta})$, we have
\begin{align}
& E\left\{ (\hat{w}_{di}^{E}-\hat{w}_{di})(\hat{w}_{dj}-w_{dj})1_{\cB%
^{c}}\right\}   \label{EexpyEy} \\
& \quad =E\left[ \left\{ \exp (\hat{y}_{di}^{E}+\hat{\alpha}_{d})-\exp (\hat{%
y}_{di}+\alpha _{d})\right\} \left\{ \exp (\hat{y}_{dj}+\alpha _{d})-\exp
(y_{dj})\right\} 1_{\cB^{c}}\right]   \notag \\
& \quad \leq E\left[ \exp (\hat{y}_{di}^{E}+\hat{y}_{dj}+\hat{\alpha}%
_{d}+\alpha _{d})1_{\cB^{c}}\right] +E\left[ \exp (\hat{y}%
_{di}+y_{dj}+\alpha _{d})1_{\cB^{c}}\right] . \notag
\end{align}%
For $\nu \in (0,1)$, we define the neighborhood $N(\btheta_{0})=\{\btheta\in
\Theta :|\btheta-\btheta_{0}|<D^{-\nu /2}\}$. Using (\ref{bdiy}) and
applying Hölder's inequality, the first expectation on the right-hand side
of (\ref{EexpyEy}) can be bounded as
\begin{align*}
& \hspace{-0.8cm}E\left[ \exp (\hat{y}_{di}^{E}+\hat{y}_{dj}+\hat{\alpha}%
_{d}+\alpha _{d})1_{\cB^{c}}\right] \leq \exp \left\{ 2\sup_{N(\btheta%
_{0})}\alpha _{d}(\btheta)\right\}  \\
& \hspace{-0.8cm}\quad \times E\left[ \exp \left\{ \sup_{N(\btheta_{0})}(\bb%
_{di}(\btheta)+\bb_{dj}(\btheta))^{\prime }\by_{s}\right\} 1_{\cB^{c}}\right]
\\
& \hspace{-0.8cm}\quad \leq \exp \left\{ 2\sup_{N(\btheta_{0})}\alpha _{d}(%
\btheta)\right\} E^{1/2}\left[ \exp \left\{ 2\sup_{N(\btheta_{0})}(\bb_{di}(%
\btheta)+\bb_{dj}(\btheta))^{\prime }\by_{s}\right\} \right] P^{1/2}(\cB%
^{c}).
\end{align*}%
But the suprema of $|\alpha _{d}(\btheta)|$ and $|\bb_{di}(\btheta)|$ over $%
N(\btheta_{0})$ are bounded. Moreover, since $\by_{s}$ is normally
distributed, the expected value on the right-hand side of the inequality is
bounded. Now by Lemma 1 of Molina (2009) with $\nu =\eta \in (0,3/4)$ and $%
b>16$, we get $P^{1/2}(\cB^{c})=O(D^{-b/16})=o(D^{-1})$.
Therefore,
\begin{equation}
E\left[ \exp (\hat{y}_{di}^{E}+\hat{y}_{dj}+\hat{\alpha}_{d}+\alpha _{d})1_{%
\cB^{c}}\right] =o(D^{-1}).  \label{EyEy2a}
\end{equation}%
Similarly, we have
\begin{align}
& E\left[ \exp (\hat{y}_{di}+y_{dj}+\alpha _{d})1_{\cB^{c}}\right]
\label{Ehatyy}\\
& \leq \exp (\alpha _{d})E^{1/2}\left[ \exp \left\{ \sup_{N(\btheta%
_{0})}\bb_{di}^{\prime}(\btheta)\by_{s}+y_{dj}\right\} \right]
P^{1/2}(\cB^{c})=o(D^{-1}).  \notag
\end{align}%
Replacing (\ref{EyEy2a}) and (\ref{Ehatyy}) in (\ref{EexpyEy}), we obtain
$E\left\{ (\hat{w}_{di}^{E}-\hat{w}_{di})(\hat{w}_{dj}-w_{dj})1_{\cB%
^{c}}\right\} =o(D^{-1})$.\\
\\
\textbf{Proof of (D):} Consider the first term on the left-hand side of (\ref{termssumM}), given by
 $$
 E\left\{\bh_{di}'\cF^{-1}\bs(\hat w_{dj}-w_{dj})\right\}=E\left(\bh_{di}'\cF^{-1}\bs\,\hat w_{dj}\right)-E\left(\bh_{di}'\cF^{-1}\bs\,w_{dj}\right)
 $$
 Using $w_{di}=\exp(\bx_{di}^{\prime}\bbeta+u_d+e_{di})$ and taking into account that
 \beq\label{hdi}
 \bh_{di}=\exp(\delta_{di})\partial\delta_{di}/\partial\btheta, \quad \delta_{di}=\alpha_d+\bx_{di}'\bbeta+\bb_{di}'\bv_s,
 \eeq
 we obtain
 \beq\label{Ehfsw}
 E\left(\bh_{di}'\cF^{-1}\bs\, \hat w_{dj}\right)=\exp\left(\alpha_d+\bx_{dij}'\bbeta\right)E\left\{\exp(\bb_{dij}^{\prime}\bv_s)\left(\partial\delta_{di}/\partial\btheta\right)^{\prime}\cF^{-1}\bs\right\}.
 \eeq
 where $\bb_{dij}=\bb_{di}+\bb_{dj}=2\eeta_d+\ff_{di}+\ff_{dj}$, with $|\eeta_d|=O(1)$ and $|\ff_{di}|=O(D^{-1/2})$ by (\ref{decompbdi}).

 To calculate the expected value in (\ref{Ehfsw}), note that $\delta_{di}=\alpha_d+\bx_{di}'\bbeta+\bb_{di}'\bv_s$ and define
 \beq\label{gC}
 \bg_d=\cF^{-1}\frac{\partial \alpha_d}{\partial\btheta}=(g_{d1},g_{d2})',\quad C_{di}=\cF^{-1}\frac{\partial \bb_{di}'}{\partial\btheta}=(\bc_{di1},\bc_{di2})'.
 \eeq
 Then, we can write
 \beq\label{Fdelta}
 \cF^{-1} \frac{\partial\delta_{di}}{\partial\btheta}=\cF^{-1}\frac{\partial \alpha_d}{\partial\btheta}+\cF^{-1}\frac{\partial \bb_{di}'}{\partial\btheta}\bv_s
 =\bg_d+C_{di}\bv_s,
 \eeq
 Moreover, denoting $\bA_h=\bP_s\bDelta_h\bP_s$, $q_h=\bv_s'\bA_h\bv_s$, $h=1,2$ and $\bq=(q_1,q_2)'$, the vector of scores (\ref{score}) can be expressed as
 \beq\label{sqnu}
 \bs=(\bq-E\bq)/2+\bnu,\quad \bnu=(\nu_1,\nu_2)',\quad \nu_h=\left\{\tr(\bP_s\bDelta_h)-\tr(\bV_s^{-1}\bDelta_h)\right\}/2.
 \eeq
 Using these expressions, we get
 \begin{eqnarray*}
 \hspace{-1 cm}&& E\left\{\exp(\bb_{dij}^{\prime}\bv_s)\left(\frac{\partial\delta_{di}}{\partial\btheta}\right)'\cF^{-1}\bs\right\}=\frac{1}{2}\bg_d^{\prime}E\left\{\exp(\bb_{dij}^{\prime}\bv_s)(\bq-E\bq)\right\}\\
 \hspace{-1 cm}&& +\bg_d^{\prime}\bnu E\left\{\exp(\bb_{dij}^{\prime}\bv_s)\right\}+\frac{1}{2}E\left\{\exp(\bb_{dij}^{\prime}\bv_s)\bv_s^{\prime}C_{di}^{\prime}(\bq-E\bq)\right\}+E\left\{\exp(\bb_{dij}^{\prime}\bv_s)\bv_s^{\prime}C_{di}^{\prime}\right\}\bnu.
 \end{eqnarray*}
 Using repeatedly Lemma 5(iv) of Molina (2009), we obtain
 \begin{equation}\label{EhFsw}
 E\left(\bh_{di}^{\prime}\cF^{-1}\bs\, \hat w_{dj}\right)=E_{dij}\left\{\tr\left(\cF^{-1}\frac{\partial\eeta_d'}{\partial\btheta}\bE_{dj}\right)+\frac{1}{2}\left(\frac{\partial\alpha_d}{\partial\btheta}
  +2\frac{\partial\eeta_d^{\prime}}{\partial\btheta}\bV_s{\eeta_{dj}}\right)^{\prime}\cF^{-1}(2\bnu+\bepsilon_{dj})\right\}.
 \end{equation}
 For the expected value $E(\bh_{di}'\cF^{-1}\bs\, w_{dj})$, note that $w_{dj}=\exp(y_{dj})$, where $y_{dj}=\bx_{dj}^{\prime}\bbeta+v_{dj}$, for $v_{dj}=u_d+e_{dj}$. However, since $j\in \bar s_d$, we cannot express $y_{dj}$ in terms of $\bv_s$ as done above. In this case, we construct an extended vector $\bv_{sj}^*=(\bv_s^{\prime},v_{dj})^{\prime}$, whose distribution is $N(\cero_2,\bV_s^*)$, for
 $$
 \bV_s^*=\left(\begin{array}{cc}
 \bV_s & \sigma_u^2\bz_d\\
 \sigma_u^2\bz_d^{\prime} & \sigma_u^2+\sigma_e^2
 \end{array}\right)
 $$
 where $\bz_d=\bZ_s\bm_d$. Defining also $\bb_{di}^*=(\bb_{di}^{\prime},1)^{\prime}$, we can express
 $$
 E\left(\bh_{di}^{\prime}\cF^{-1}\bs\, w_{dj}\right)=\exp\left(\alpha_d+\bx_{dij}^{\prime}\bbeta\right)E\left[\exp\{(\bb_{di}^*)^{\prime}\bv_{sj}^*\}\left(\frac{\partial\delta_{di}}{\partial\btheta}\right)^{\prime}\cF^{-1}\bs\right].
 $$
 Expressing now $\cF^{-1} \partial\delta_{di}/\partial\btheta$ and $\bs$ in terms of $\bv_{sj}^*$ similarly as in (\ref{Fdelta}) and (\ref{sqnu}) by adding zero elements to the vectors and matrices multiplying $\bv_s$, we can apply exactly the same results as used for $E\left(\bh_{di}'\cF^{-1}\bs\, \hat w_{dj}\right)$. The result turns out to be equal to (\ref{EhFsw}) with $E_{dij}$ replaced by $E_{dij}^*$.

 The rest of terms on the left-hand side of (\ref{termssumM}) are obtained following a similar procedure, by expressing the terms within the expectations as sums of products of quadratic and linear forms in $\bv_s$ multiplied by exponentials of linear forms of $\bv_s$ and then applying repeatedly Lemma 5 of Molina (2009).
\vspace{0.3cm}\newline
\textbf{Proof of (E):} Note that
 \begin{equation}\label{difEabs}
 E\left\{\left|\bh_{di}'\cF^{-1}\bs(\hat w_{dj}-w_{dj})\right|1_{\cB^c}\right\}\leq E\left\{\left|\bh_{di}'\cF^{-1}\bs\,\hat w_{dj}\right|1_{\cB^c}\right\}
 +E\left\{\left|\bh_{di}'\cF^{-1}\bs\, w_{dj}\right|1_{\cB^c}\right\}.
 \end{equation}
 By the definition of $\bh_{di}$ in (\ref{hdi}) and that of $\hat w_{dj}$ in (\ref{what}), we obtain
 \begin{eqnarray*}
 E\left\{\left|\bh_{di}'\cF^{-1}\bs\,\hat w_{dj}\right|1_{\cB^c}\right\}=\exp(2\alpha_d+\bx_{dij}'\bbeta)E\left[\exp\{(\bb_{di}+\bb_{dj})'\bv_s\}\left|\left(\frac{\partial \delta_{di}}{\partial\btheta}\right)'\cF^{-1}\bs\right|1_{\cB^c}\right].
 \end{eqnarray*}
 Now applying repeatedly H\"older's inequality, we get
 \begin{eqnarray}
 && E\left\{\left|\bh_{di}'\cF^{-1}\bs\,\hat w_{dj}\right|1_{\cB^c}\right\}\leq \exp(2\alpha_d+\bx_{dij}'\bbeta)E^{1/2}\left[\exp\{2(\bb_{di}+\bb_{dj})'\bv_s\}\right]\nonumber\\
 && \quad\times E^{1/8}\left|\frac{\partial \delta_{di}}{\partial\btheta}\right|^8E^{1/8}\left|\cF^{-1}\bs\right|^8P^{1/4}(\cB^c)=O(D^{-1/2-b/32})=o(D^{-1})\label{EhFshatw1}
 \end{eqnarray}
 for $b>16$, noting that by the proof of Theorem 1 in Molina (2009), it holds
 \begin{equation}\label{E8dE8Fs}
 E^{1/8}\left|\frac{\partial \delta_{di}}{\partial\btheta}\right|^8=O(1),\quad E^{1/8}\left|\cF^{-1}\bs\right|^8=O(D^{-1/2}),
 \end{equation}
 that $P^{1/4}(\cB^c)=O(D^{-b/32})$, by Lemma 1 in Molina (2009) with $\nu=\eta\in (0,3/4)$, and finally taking into account that $\bv_s$ is normally distributed and that
 $\exp(2\alpha_d+\bx_{dij}'\bbeta)$ and $\bb_{di}$ are bounded.
 By a similar reasoning, we obtain
 \begin{equation}\label{EhcFw1}
 E\left\{\left|\bh_{di}'\cF^{-1}\bs\, w_{dj}\right|1_{\cB^c}\right\}\leq \exp(\alpha_d+\bx_{dij}'\bbeta)
E\left[\exp\{(\bb_{di}^*)^{\prime}\bv_{sj}^*\}\left|\left(\frac{\partial \delta_{di}}{\partial\btheta}\right)'\cF^{-1}\bs\right|1_{\cB^c}\right]=o(D^{-1}).
 \end{equation}
 By (\ref{EhcFw1}) and (\ref{EhFshatw1}), we obtain
 $E\left\{\left|\bh_{di}'\cF^{-1}\bs(\hat w_{dj}-w_{dj})\right|1_{\cB^c}\right\}=o(D^{-1})$.
 The remaining results in (E) are proved similarly.
\hfill $\Box $



\subsection*{PROOF OF THEOREM \protect\ref{them7}}

Similarly as before, we spell the proof for ML, since for REML estimation the proof is analogous.
For $\nu \in (0,1)$, let us define the neighborhood
\begin{equation*}
N(\bbeta_{0},\btheta_{0})=\{(\bbeta^{\prime },\btheta^{\prime })^{\prime
}\in \Theta \times \RR^{p};|\bbeta-\bbeta_{0}|<D^{-\nu /2},\ |\btheta-\btheta%
_{0}|<D^{-\nu /2}\}.
\end{equation*}%
By a first-order Taylor expansion of $M_{3d,ij}(\bbeta,\btheta)$ around $(%
\bbeta,\btheta)=(\bbeta_{0},\btheta_{0})$ evaluated at the ML estimates $(%
\hat{\bbeta},\hat{\btheta})$, we obtain
\begin{equation}
M_{3d,ij}(\hat{\bbeta},\hat{\btheta}) =M_{3d,ij}(\bbeta_{0},\btheta%
_{0})+\left. \frac{\partial M_{3d,ij}(\bbeta,\btheta)}{\partial \btheta}%
\right\vert _{(\bbeta_{\ast },\btheta_{\ast })}(\hat{\btheta}-\btheta_{0}) +\left. \frac{\partial M_{3d,ij}(\bbeta,\btheta)}{\partial \bbeta}%
\right\vert _{(\bbeta_{\ast },\btheta_{\ast })}(\hat{\bbeta}-\bbeta_{0}),
\end{equation}%
where $(\bbeta_{\ast }^{\prime },\btheta_{\ast }^{\prime })^{\prime }\in N(%
\bbeta_{0},\btheta_{0})$. Taking expected value, we obtain
\begin{align}
& \hspace{-1cm}E\left[ M_{3d,ij}(\hat{\bbeta},\hat{\btheta})\right]
=M_{3d,ij}(\bbeta_{0},\btheta_{0})  \label{EM3dij} \\
& \hspace{-1cm}+E\left[ \left. \frac{\partial M_{3d,ij}(\bbeta,\btheta)}{%
\partial \btheta}\right\vert _{(\bbeta_{\ast },\btheta_{\ast })}(\hat{\btheta%
}-\btheta_{0})+\left. \frac{\partial M_{3d,ij}(\bbeta,\btheta)}{\partial %
\bbeta}\right\vert _{(\bbeta_{\ast },\btheta_{\ast })}(\hat{\bbeta}-\bbeta%
_{0})\right].  \notag
\end{align}%
where we have
\begin{align}
& \hspace{-1cm}E\left[ \left. \frac{\partial M_{3d,ij}(\bbeta,\btheta)}{%
\partial \btheta}\right\vert _{(\bbeta_{\ast },\btheta_{\ast })}(\hat{\btheta%
}-\btheta_{0})+\left. \frac{\partial M_{3d,ij}(\bbeta,\btheta)}{\partial %
\bbeta}\right\vert _{(\bbeta_{\ast },\btheta_{\ast })}(\hat{\bbeta}-\bbeta%
_{0})\right]  \label{EM3ij} \\
& \hspace{-1cm}\leq \left( \sup_{N(\bbeta_{0},\btheta_{0})}\left\vert \frac{%
\partial M_{3d,ij}(\bbeta,\btheta)}{\partial \btheta}\right\vert \right) E|%
\hat{\btheta}-\btheta_{0}|+\left( \sup_{N(\bbeta_{0},\btheta_{0})}\left\vert
\frac{\partial M_{3d,ij}(\bbeta,\btheta)}{\partial \bbeta}\right\vert
\right) E|\hat{\bbeta}-\bbeta_{0}|.  \notag
\end{align}%
By Lemma 1 in Molina (2009), for every $\nu \in (0,1)$, we have
$\hat{\btheta}-\btheta_{0}=\cF^{-1}\bs+\br^{\ast }$,
where $|\br^{\ast }|\leq D^{-\nu }E(u^{\ast })$, where $E(u^{\ast })=O(1)$;
hence, $|\br^{\ast }|=O(D^{-\nu })$. As a consequence, we have
\begin{equation*}
E|\hat{\btheta}-\btheta_{0}|\leq E|\cF^{-1}\bs|+E|\br^{\ast }|,
\end{equation*}%
and since $E|\cF^{-1}\bs|=O(D^{-1/2})$, we obtain that
\begin{equation}
E|\hat{\btheta}-\btheta_{0}|=O(D^{-1/2-\nu }),\ \nu \in (0,1).
\label{Ehattheta}
\end{equation}%
Note also that
$\hat{\bbeta}-\bbeta_{0}=\bQ_{s}(\hat{\btheta})\bX_{s}^{\prime }\bV_{s}^{-1}(%
\hat{\btheta})\bv_{s}$.
Then, we can write
\begin{equation*}
E|\hat{\bbeta}-\bbeta_{0}|=\left( \sup_{N(\bbeta_{0},\btheta_{0})}\Vert \bQ%
_{s}(\btheta)\Vert \right) \left( \sup_{N(\bbeta_{0},\btheta_{0})}\Vert \bV%
_{s}^{-1}(\btheta)\Vert \right) \Vert \bX_{s}\Vert E|\bv_{s}|.
\end{equation*}%
By Lemma \ref{lemmaV} (ii) and (iii), we know that at the true value of $%
\btheta$, $\Vert \bQ_{s}\Vert =O(D^{-1})$ and $\Vert \bV_{s}^{-1}\Vert =O(1)$%
. By continuity of $\bQ_{s}(\btheta)$ and $\bV_{s}^{-1}(\btheta)$ on $%
\btheta
$, we have
\begin{equation*}
\sup_{N(\btheta_{0},\bbeta_{0})}\Vert \bQ_{s}(\btheta)\Vert =O(D^{-1}),\quad
\sup_{N(\btheta_{0},\bbeta_{0})}\Vert \bV_{s}^{-1}(\btheta)\Vert =O(1).
\end{equation*}%
Considering the facts that $\Vert \bX_{s}\Vert =O(D^{1/2})$ and $E|\bv%
_{s}|=O(1)$, we obtain
\begin{equation}
E|\hat{\bbeta}-\bbeta_{0}|=O(D^{-1/2}).  \label{Ehatbeta}
\end{equation}%
By replacing (\ref{Ehatbeta}) and (\ref{Ehattheta}) in (\ref{EM3ij}), the
desired result is obtained if the following conditions hold:
\begin{equation*}
\sup_{N(\bbeta_{0},\btheta_{0})}\left\vert \frac{\partial M_{3d,ij}(\bbeta,%
\btheta)}{\partial \btheta}\right\vert =O(D^{-1/2}),\ \sup_{N(\bbeta_{0},%
\btheta_{0})}\left\vert \frac{\partial M_{3d,ij}(\bbeta,\btheta)}{\partial %
\bbeta}\right\vert =o(D^{-1/2}).
\end{equation*}%
Now write (\ref{M3dij}) as
$M_{3d,ij}(\bbeta,\btheta)=\left\{E_{dij}(\bbeta,\btheta)-E_{dij}^{\ast }(\bbeta,\btheta)\right\}M_{31,ij}(\btheta)$,
where $M_{31,ij}(\btheta)=K_d(\btheta)/2+C_d(\btheta)$.
Now since $M_{31,ij}(\btheta)$ does not depend on
$\bbeta$ and
\begin{equation*}
\frac{\partial E_{dij}(\bbeta,\btheta)}{\partial \bbeta}=E_{dij}(\bbeta,%
\btheta)\bx_{dij},\quad \frac{\partial E_{dij}^{\ast }(\bbeta,\btheta)}{%
\partial \bbeta}=E_{dij}^{\ast }(\bbeta,\btheta)\bx_{dij},
\end{equation*}%
Then, we have
\begin{equation*}
\left\vert \frac{\partial M_{3d,ij}(\bbeta,\btheta)}{\partial \bbeta}%
\right\vert \leq \left\{E_{dij}(\bbeta,\btheta)+E_{dij}^{\ast }(\bbeta,\btheta)\right\}|\bx_{dij}||M_{31,ij}(\btheta%
)|.
\end{equation*}%
Therefore,
\begin{equation}
\sup_{N(\bbeta_{0},\btheta_{0})}\left\vert \frac{\partial
M_{3d,ij}(\bbeta,\btheta)}{\partial \bbeta}\right\vert  \leq |\bx_{dij}|\sup_{N(%
\bbeta_{0},\btheta_{0})}\left\{E_{dij}(\bbeta,\btheta)+E_{dij}^{\ast }(\bbeta,\btheta)\right\}
\sup_{N(\bbeta_{0},\btheta_{0})}|M_{31,ij}(\btheta)|.\label{supM3dij}
\end{equation}%
We know that $|\bx_{dij}|=O(1)$. Moreover, it is easy to see that the
suprema over $N(\btheta_{0},\bbeta_{0})$ of $E_{dij}(\bbeta,\btheta)+
E_{dij}^{\ast }(\bbeta,\btheta)$ is bounded. Finally, it is also easy but
cumbersome to check that
\begin{equation*}
\sup_{N(\btheta_{0},\bbeta_{0})}|M_{31,ij}(\btheta)|=O(D^{-1}).
\end{equation*}%
By (\ref{supM3dij}), this implies
\begin{equation*}
\sup_{N(\btheta_{0},\bbeta_{0})}\left\vert \frac{\partial M_{3d,ij}(\bbeta,%
\btheta)}{\partial \bbeta}\right\vert =O(D^{-1}).
\end{equation*}%
It also holds that
\begin{equation}
\sup_{N(\btheta_{0},\bbeta_{0})}\frac{\partial E_{dij}(\bbeta,\btheta)}{%
\partial \btheta}=O(1),\quad \sup_{N(\btheta_{0},\bbeta_{0})}\frac{\partial
E_{dij}^{\ast }(\bbeta,\btheta)}{\partial \btheta}=O(1)  \label{supEdij}
\end{equation}%
and that
\begin{equation}
\sup_{N(\btheta_{0},\bbeta_{0})}\left\vert \frac{\partial M_{31,ij}(\bbeta,%
\btheta)}{\partial \btheta}\right\vert =O(D^{-1/2}).  \label{supMdij}
\end{equation}%
Relations (\ref{supEdij}) and (\ref{supMdij}) imply that
\begin{equation*}
\sup_{N(\bbeta_{0},\btheta_{0})}\left\vert \frac{\partial M_{3d,ij}(\bbeta,%
\btheta)}{\partial \btheta}\right\vert =O(D^{-1/2}).
\end{equation*}%
Finally, (\ref{EM3dij}) and (\ref{EM3ij}) lead to
\begin{equation*}
E\left[ M_{3d,ij}(\hat{\bbeta},\hat{\btheta})\right] =M_{3d,ij}(\bbeta_{0},%
\btheta_{0})+o(D^{-1}),
\end{equation*}%
which is our desired result. \hfill $\Box $\newline

\end{document}